\documentclass{article}
\usepackage[utf8]{inputenc} % utf8 coding
\usepackage[a4paper, total={6in, 8in}]{geometry}
\usepackage{amssymb,amsmath} % math related stuff
\usepackage{graphicx} % including graphics
\usepackage{authblk}

\usepackage{array} % array package
\usepackage{eurosym} % the euro-symbol
\usepackage{xcolor} % colored text
\usepackage{booktabs}
\usepackage{epstopdf}

\usepackage{algorithm} 
\usepackage{algpseudocode}
\usepackage{multicol}
\usepackage{rotating}
\usepackage{multirow}
\usepackage{subcaption}
\usepackage{natbib}
\usepackage{hyperref}
\usepackage{adjustbox}
\usepackage{graphicx}
\usepackage{multirow}
\usepackage{amsmath} 
\usepackage{setspace}

\begin{document}

\title{A generic stochastic network flow formulation for production optimization of district heating systems}

\author[1,2,*]{Daniela Guericke}
\author[2]{Amos Schledorn}
\author[2]{Henrik Madsen}
\affil[1]{\footnotesize High-Tech Business and Entrepreneurship Department, University of Twente, P.O. Box 217, 7500 AE, Enschede, The Netherlands}
\affil[2]{Technical University of Denmark, Department of Applied Mathematics and Computer Science, Richard Petersens Plads, 2800 Kgs. Lyngby, Denmark}
\affil[*]{Corresponding author, \textit {d.guericke@utwente.nl}}
%\eid{123@gmail.com}

\maketitle
%\thanks{Ignacio Blanco and Daniela Guericke (corresponding author)are with the Technical University of Denmark, DK-2800 Kgs. Lyngby, Denmark (email addresses: \{igbl, dngk\}@dtu.dk). Anders N. Andersen is with EMD International A/S (email: ana@emd.dk).  The work of all authors is partly funded by DSF (Det Strategiske Forskningsr{\aa}d) through the CITIES research center (no. 1035-00027B).}}

\begin{abstract}
	
	District heating is an important component in the EU strategy to reach the set emission goals, since it allows an efficient supply of heat while using the advantages of sector coupling between different energy carriers such as power, heat, gas and biomass. Most district heating systems use several different types of units to produce heat for hundreds or thousands of households. The technologies reach from natural gas-fired and electric boilers to biomass-fired units as well as waste heat from industrial processes and solar thermal units. Furthermore, combined heat and power units (CHP) units are often included to use the synergy effects of excess heat from electricity production. 
	We propose a generic mathematical formulation for the operational production optimization in district heating systems. The generality of the model allows it to be used for most district heating systems although they might use different combinations of technologies in different system layouts. The mathematical formulation is based on stochastic programming to account for the uncertainty of production from non-dispatchable units such as waste heat and solar heat. Furthermore, the model is easily adaptable to different application cases in district heating such as operational planning, bidding to electricity markets and long-term evaluation. We present results from three real cases in Denmark with different requirements.
	
\end{abstract}

\section{Introduction} \label{sec:introduction}
Under the European Green Deal, the European Commission aims at net carbon neutrality by 2050 \citep{europeancommission2019}, which requires the transformation of the European energy system through integration of a large amount of renewable energy \citep{Hainsch2022}. District heating (DH) can play a substantial role in supporting this transition. Not only is heating and cooling responsible for half of the EU's final energy consumption \citep{europeancommission2015}, the flexibility potential of DH systems can also ease the integration of renewable electricity sources substantially \citep{Thomassen2021}. The Danish government specifically aims at 70 \% emission reduction by 2030 (compared to 1990) alongside carbon neutrality by the middle of the century \citep{Folketinget2020}. This transition requires smart energy systems and modelling approaches, that do not merely focus on a single sector, but take different energy carriers into account \citep{Lund2017a}. This modelling paradigm is particularly important for the practical ability of DH optimization models, as systems often feature a variety of energy sources, such as biomass, solar and industrial waste heat \citep{stateofgreen2018}, and the optimal operation of some systems must even take behind-the-meter electricity dispatch into account \citep{hvidesande}.

Following this notion, this article proposes a novel generic optimisation model for DH systems. Our model features a network flow formulation based on stochastic programming that can take wide variety of energy carriers, productions units, markets and demand sinks into account. Furthermore, the model can be used in different application cases such as operational planning, bidding and long-term system analysis by merely changing input data and the non-anticipativity constraints of the stochastic model. The applicability of the model to all three cases is shown based on real data from the three DH systems in Denmark.

The remainder of this paper is structured as follows. Related work and our contributions are presented in Section \ref{sec:literature}. Sections \ref{sec:network} and \ref{sec:model} present the network and mathematical formulation of our proposed optimization model, respectively. The DH systems we use as cases are introduced in Section \ref{sec:cases} and numerical results are presented in Section \ref{sec:results}. Finally, Section \ref{sec:conclusion} summarizes our work and gives an outlook on future research.

\section{Related work} \label{sec:literature}
We propose a network-flow based model for the optimal scheduling of different energy generation and conversion units in DH systems. Our formulation is based on stochastic programming and allows to model an arbitrary range of energy carriers. Hence, related work focuses on mathematical models that (1) can optimise operational scheduling in DH networks (2) under uncertainty and (3) can represent arbitrary energy carriers. To the best of our knowledge, existing research meeting all of these criteria is limited. In the following, we provide a brief review of other studies focusing on energy system models (Section \ref{subsec:lit-esoms}), DH models (Section \ref{subsec:lit-dhModels}) and energy hubs (Section \ref{subsec:lit-energyHubs}).

\subsection{Energy System Modelling Frameworks} \label{subsec:lit-esoms}
Energy system optimisation models (ESOMs) are most commonly formulated as large-scale linear programs with the aim of providing the optimal dispatch of one or more energy carriers \citep{Dominkovic2018f} and/or investments \citep{Dominkovic2020} in related technologies. The scope and scale of applications typically reach from the district \citep{Weckesser2021}, urban \citep{Dominkovic2020} to country \citep{Daly2014} or even continental level \citep{Pavicevic2020a}, spanning one or more years. An overview on different modelling frameworks can be found on the OpenMod Wiki \citep{Openmod}. A focus on local multi-energy systems is provided in \citep{cuisinier2021techno}.
The model proposed in this paper is a stochastic program, taking uncertainty modelled as scenarios into account, and is formulated in a generic way such that it is able to model arbitrary energy commodities and technologies. Notably, most energy system modelling frameworks are either  deterministic models (e.g. Balmorel \citep{Wiese2018a}) or model certain energy technologies specifically \cite{Helisto2019} thus not reaching the general applicability this work does. The Balmorel extension OptiFlow by \cite{Ravn2017} does, as the model proposed here, use a graph-based network flow formulation, but is, at the time of writing, focused on the waste sector. \citep{hilpert2018open} present the oemof modelling framework that uses a network flow formulation to model an energy system, similar to the setting used in this paper, but without consideration of uncertainty and market interaction. \citep{Quelhas2007} provide a general network flow formulation for US power network incorporating also gas and coal inflow. The resulting model is a deterministic LP for long-term evaluation that abstracts from operational constraints such as on/off status, market interaction and renewable energy sources (RES). Also other frameworks, such as Calliope by \cite{Pfenninger2018} meet comparable requirements in theory, but are mostly applied to long-term planning problems (e.g. \cite{Pickering2019}) rather than having an operational focus. That focus makes them less suitable for real-world operational problems, where features as uncertainty modelling, rolling horizon optimisation and market interaction are important.

\subsection{District Heating Optimisation Models}\label{subsec:lit-dhModels}
District-heating specific models are often able to capture detailed characteristics of the system and  take uncertainty into account: In \cite{Zhou2020a}, a distributionally robust linear formulation of a co-generation dispatch problem for heat and power dispatch is proposed. The authors of \cite{Xue2021} solve a robust unit commitment problem for a co-generating heat and power system. In \cite{Hohmann2019}, operation strategies for a DH system are optimised in a stochastic program. However, such models, coming along with rich detail with respect to district-heating specific characteristics, typically model energy carriers and units, usually heat and possibly power, explicitly, thus lacking general applicability to systems with a wider range of energy carriers.
 
Co-generation of power and heat raises the question of an integrated optimisation of both dispatch and power market bidding, both on day-ahead \cite{hurb} and regulating power markets \citep{hvidesande} and taking different bidding products into account \citep{blockbids}. Optimisation under uncertainty becomes especially useful in such co-generation problems, where typical sources of uncertainty can include heat load \citep{Dimoulkas2015}, prices \citep{blockbids} or RES generation \citep{hvidesande}. These authors model the studied DH systems for bidding but abstain from a generic model formulation. An example of a network flow formulation in DH is the work in \cite{Bordin2016}, where the planning of an Italian DH system is optimised, while DH characteristics are modelled explicitly but no uncertainty is considered. The reader is referred to \cite{Sarbu2019} for an review of DH system optimisation models.

Operational optimisation in DH is often combined with design or capacity expansion models. Examples are \cite{Wirtz2020}, \cite{gabrielli2018optimal} and \cite{Weinand2019} relying on explicit modelling of the systems. Further publications including investment decisions in multi-energy systems are presented in \citep{cuisinier2021techno}.

\subsection{Energy hubs}\label{subsec:lit-energyHubs}
The energy hub concept was defined by \cite{Geidl2007} as "An energy hub is considered a unit where multiple energy carriers can be converted, conditioned, and stored", i.e, an energy hub can contain several energy units, storage facilities and transformers. Several energy hubs are  interconnected using transmission systems. For each energy hub, a so-called coupling matrix can be derived, which is used to determine the operational strategy\citep{Geidl2007,Beigvand2017}. This matrix contains the transformation factors from one energy carrier to another energy carrier based on the setup of all units inside the energy hub. Based on the energy hub concept, several planning tasks can be executed, e.g., design of an energy hub  \cite{wang2018mixed}, systems impact analysis of energy hubs, operational planning  \cite{Geidl2005, Najafi2016} and optimal power flow. Within the context of this publication, the productions units of the DH operator can be considered as one energy hub. We are interested in a generic formulation for the stochastic operational scheduling inside the energy hub. We also use a transformation matrix similar to the coupling matrix, but on a unit-level to define the transformations of each unit. \cite{Beigvand2017} look at the economic dispatch of energy hub and give many examples for possible energy hubs and a mathematical formulation using the coupling matrix. Their formulation abstracts from specific time periods, commitment decisions and uncertainties, since it is concerned with optimizing the energy input and dispatch factors to reach the required output. \cite{Geidl2005} propose a general formulation for the determining the optimal dispatch of the units inside an energy hub for one hour using a non-linear formulation. \cite{Najafi2016} present a specific stochastic model for an energy hub containing CHP units, generators and wind turbines considering electricity prices and wind power production as uncertain. \cite{wang2018mixed} determine the configuration of an energy hub using a model that uses a general notation based on input and output ports and branches between components. The decisions are the possible combinations as binary variables and the optimal operation strategy.  An extended overview of models for energy hubs is given in \cite{Mohammadi2017}. Like in the above mentioned publications, most energy hubs publication focus on planning and analysis and do not consider real-time operational planning \cite{Krause}.

\subsection{Contribution}
Based on the presented literature above, we summarize our contributions as follows:
\begin{enumerate}
    \item We propose a novel formulation for the optimisation of DH systems. The model relies on a network representation of the DH system, which makes the optimization model itself very generic and it can be easily applied to a wide range of DH systems with different units and requirements. The model uses stochastic programming to account for uncertainty in prices and heat flows. 
    \item The proposed model can be used in different application cases. Traditional operational optimization dispatching the heat production units can achieved by using the model with non-anticipativity on the commitment status and production on some of the units. Determination of bids to day-ahead markets can be achieved by including non-anticipativity constraints that create curves based on the electricity prices scenarios (based on the work in \cite{pandvzic2013offering, hvidesande}). Finally, a deterministic version of the model can be used for pure evaluation purposes on historic data.  
    \item The optimisation model also allows defining sliding time windows and using the model in a rolling or receding horizon approach.
    \item The performance of the model in terms of costs, energy mixes and runtimes is extensively evaluated in several case studies using real data from three Danish DH networks including out-of-sample testing.
\end{enumerate}

\section{Network representation} \label{sec:network}
In this section, we present how a DH system can be represented as a network graph with arcs and vertices and introduce all components, parameters and sets. The general idea is to transform all components of a DH system such as production units, storage units, demand sites etc. to vertices in a network and use the arcs to model possible flows of energy between units. An overview of the nomenclature is given in Table \ref{tab:nomenclature_dh}. The model formulation in Section \ref{sec:model} uses this structure as basis for the mathematical formulation. 
\begin{table}[ht]
\footnotesize
    \caption{Nomenclature - DH network components}
    \centering
    \begin{adjustbox}{width=\textwidth}
    \begin{tabular}{p{0.15\columnwidth}p{0.85\columnwidth}}\toprule
    Symbol & Description\\\midrule
         $\mathcal{U}$&  Production units\\
         $\mathcal{E}$& Energy sources \\
         $\mathcal{D}$ & Demand sites\\
         $\mathcal{P}$ & Pipeline connection\\
         $\mathcal{S}$ & Storage units\\
         $\mathcal{F}$ & Energy type\\
           $\mathcal{T}$ & Periods\\
          ${\Omega}$ & Scenarios\\
         $\mathcal{V}$ & Vertices $\mathcal{V} = \mathcal{U} \cup \mathcal{E} \cup \mathcal{D} \cup \mathcal{P} \cup \mathcal{S}$\\
		$\mathcal{U}^{\text{WC}}$ & Set of unit vertices needing commitment decisions, $\mathcal{U}^{\text{WC}}\subseteq{\mathcal{U}}$\\
        $\mathcal{U}^{\text{*}}$ & Set of units with here-and-now decisions $\mathcal{U}^{\text{*}} \subseteq{\mathcal{U}}$\\
        $\mathcal{U}^{EXC}_u$ & Set of unit vertices excluded from production, if unit $v$ is producing, $\mathcal{U}^{EXC}_u\subseteq{\mathcal{U}}$\\
         $\mathcal{U}^{DEP}_u$ & Set of unit vertices also producing, if unit $v$ is producing, $\mathcal{U}^{DEP}_u\subseteq{\mathcal{U}}$\\
		$\mathcal{A}$ & Set of arcs, $\mathcal{A} \subset (V \times V \times F \times T \times T \times \Omega)$\\
		$\mathcal{A}^{OUT}_{v,f,t,\omega}$ & Set of outgoing arcs with energy type $f \in \mathcal{F}$ from vertex $v \in \mathcal{V}$ in period $t \in \mathcal{T}$ and scenario $\omega \in \omega$\\
		$\mathcal{A}^{IN}_{v,f,t,\omega}$ & Set of incoming arcs with energy type $f \in \mathcal{F}$ from vertex $v \in \mathcal{V}$ in period $t \in \mathcal{T}$ and scenario $\omega \in \omega$\\
		$\mathcal{M} = \mathcal{M}^\text{B} \cup \mathcal{M}^\text{S}$ & Set of buying (B) and selling (S) markets\\\midrule
         $\phi_{c,f,f'}$ & Conversion factor of vertex $v \in \mathcal{V}$ from energy type $f \in \mathcal{F}$ to energy type $f' \in \mathcal{F}$\\
         $\underline{I}_{v,f,t,\omega}/\overline{I}_{v,f,t,\omega}$ & Lower/upper bound on input of energy type $f \in \mathcal{F}$ to vertex $v \in \mathcal{V}$ in period $t \in \mathcal{T}$ and scenario $\omega \in {\Omega}$\\
           $L_{v,v',f}$ & Upper bound on flow of energy type $f \in \mathcal{F}$ from vertex $v \in \mathcal{V}$ to vertex $v' \in \mathcal{V}$\\
             $\underline{O}_{v,f,t,\omega}/\overline{O}_{v,f,t,\omega}$ & Lower, upper bound on output of energy type $f \in \mathcal{F}$ to vertex $v \in \mathcal{V}$ in period $t \in \mathcal{T}$ and scenario $\omega \in {\Omega}$\\
          $C^{I}_{v, f, t, \omega}/C^{O}_{v, f, t, \omega}$ & Cost for inflow/outflow of one unit of energy type $f \in \mathcal{F}$ at vertex $v \in \mathcal{V}$ in period $t \in \mathcal{T}$ and scenario $\omega \in \Omega$\\
            $C^{S}_{u}$ & Start-up cost for unit $u \in \mathcal{U}^{WC}$ \\
          	$T^{UT}_u/ T^{DT}_u$ & Minimum up time/down time for unit $u \in \mathcal{U}^{WC}$\\
          		$B_{u}$ & Initial online status of unit $u \in \mathcal{U}^\text{WC}$\\
		$T^\text{B}_{u}$ & Minimum remaining periods of initial online status of unit $u \in \mathcal{U}^\text{WC}$\\
		$R^U_{u,f}$, $R^D_{u,f}$ & Ramping limits on energy type $f$ for unit vertices $u \in \mathcal{U}$\\
			$c_{a}$ & Cost per unit on arc $a \in \mathcal{A}$, weighted with scenario probability\\
		$\phi_{v,f,f'}$ & Conversion factor from energy type $f$ to type $f'$ at vertex $v$\\
		$l_{a}/u_{a}$ & Lower/upper bound on flow on arc $a \in \mathcal{A}$\\
		$p_{m, t, \omega}$ & Price at market $m \in \mathcal{M}$ in period $t \in \mathcal{T}$ and scenario $\omega \in \Omega$\\\bottomrule
    \end{tabular}
    \end{adjustbox}
    \label{tab:nomenclature_dh}
\end{table}

\subsection{General sets and parameters}
Energy types are defined by the set $\mathcal{F}$ and are any kind of input fuel or output product that is used or produced in the DH network. Typical examples are electricity and heat as output products as well as electricity, wood chips, natural gas, waste heat and solar heat as input energy types.

The planning horizon is denoted by the set of periods $\mathcal{T}$. The subset of periods $\mathcal{T}^* = \{ 1, ..., |\mathcal{T}^*|\} \subseteq \mathcal{T}$ are the periods for which non-anticipativity must hold for the later defined units. Uncertain input data is given by the set of scenarios $\Omega$. Each scenario $\omega$ has a probability $\pi_\omega$ with $\sum_{\omega \in \Omega} \pi_\omega = 1$.

\subsection{Network structure}
In our model formulation, the DH system is represented by the set of vertices $\mathcal{V}$ and the set of arcs $\mathcal{A}$ that connect the vertices. Thus, we can formulate the main part of the optimization model as a flow problem on this network structure. An arc $a$ is defined by the indices $a = (v,v',f,t,t',\omega)$ where $v$ is the start vertex, $v'$ is the end vertex, $f$ is the type of energy flowing on this arc, $t$ is the start time period, $t'$ is the end time period and $\omega$ the scenario. 

To incorporate costs and restrictions on the flow in the network, vertices and arcs have several parameters. Each arc $a$ has cost per unit of flow denoted by $c_a$, which is weighted with probability $\pi_\omega$ with $\omega$ being the scenario of this arc. The flow on each arc $a$ is limited by the lower and upper bound $l_a$ and $u_a$, respectively. 

\begin{figure}[t]
    \centering
    \includegraphics[width=0.7\columnwidth]{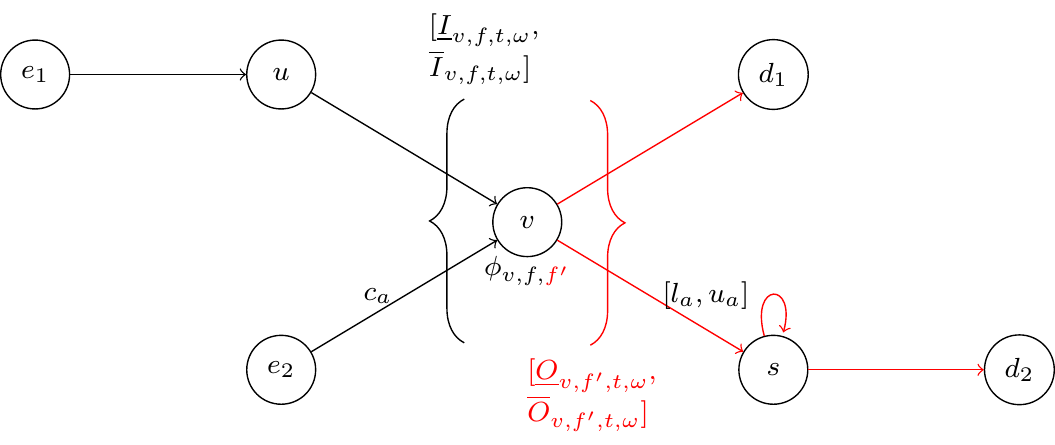}
    \caption{Example for arc and vertex parameters where vertex $v$ receives input of energy type $f$ from unit $u$ or energy source $e_2$ and transforms it to energy type $f'$ that is stored in storage unit $s$ or used in demand site $d_1$. Unit $u$ uses input from energy source $e_1$ and demand site $d_2$ gets input from storage $s$.}
    \label{fig:vertex}
\end{figure}

Each vertex $v$ has lower and upper bounds on the total incoming flow, $[\underline{I}_{v,f,t,\omega},\overline{I}_{v,f,t,\omega}]$, as well as the total outgoing flow, $[\underline{O}_{v,f,t,\omega},\overline{O}_{v,f,t,\omega}]$, per energy type $f$, period $t$ and scenario $\omega$.  Each vertex has the capability of transforming an energy type $f$ to another energy type $f'$ at a given conversion rate denoted by the transformation factor $\phi_{v,f,f'}$. The parameters are illustrated in Figure \ref{fig:vertex}.

The set of incoming and outgoing arcs for vertex $v \in \mathcal{V}$ for energy type $f \in \mathcal{F}$ in period $t\in\mathcal{T}$ and scenario $\omega \in \Omega$ are denoted by the arc subsets $\mathcal{A}^{IN}_{v,f,t,\omega}$ and $\mathcal{A}^{OUT}_{v,f,t,\omega}$, respectively.

\subsection{District heating network components} \label{sec:components}
In the following, we describe the components of the DH system and how they can be expressed using the above network structure.

\subsubsection{Energy sources $\mathcal{E}$}
Energy sources are given in the set $\mathcal{E} \subset \mathcal{V}$ and are the only possibility to insert energy of different types into the network without producing it by units. Energy sources are used for input fuels such as natural gas, biomass or electricity. Additionally, they can represent heat that is not produced but injected through waste heat sites or solar thermal units. The output of the energy source vertex $v$ is limited by  $[\underline{O}_{v,f,t,\omega},\overline{O}_{v,f,t,\omega}]$ for the energy type $f$ associated with this source and $[0,0]$ for all other energy types. The limits can vary over time and per scenario to model time-varying and/or uncertain inflow from e.g., waste heat or solar units. Furthermore, the energy can be priced with $C^O_{v,f,t,\omega}$, potentially depending on time and scenario. Thus, each arc $a$ leaving the source $v$ contains the cost $C^O_{v,f,t,\omega}$ in the parameter $c_a$.  There is no inflow to energy sources, i.e, $[\underline{I}_{v,f,t,\omega},\overline{I}_{v,f,t,\omega}] = [0,0]$. An energy source can also be used to model penalty costs for missing heat by providing heat at a high cost.

\subsubsection{Demand sites $\mathcal{D}$}
Demand sites are defined by the set $\mathcal{D} \subset \mathcal{V}$. Typical demand sites are the heat load demands in the DH network but also the electricity markets. Demand sites $v \in \mathcal{D}$ have a limitation on inflow  $[\underline{I}_{v,f,t,\omega},\overline{I}_{v,f,t,\omega}]$ for the energy type of the demand site.  The inflow for all other energy types is $[0,0]$. Demand sites have no outflow, i.e., $[\underline{O}_{v,f,t,\omega},\overline{O}_{v,f,t,\omega}] = [0,0]$. There are two common cases: $\underline{I}_{v,f,t,\omega} = \overline{I}_{v,f,t,\omega}$, if the demand needs to be fulfilled exactly, and $\underline{I}_{v,f,t,\omega} < \overline{I}_{v,f,t,\omega} = \infty$, if the demand needs to be covered but may be exceeded by its supply. The lower limit can be set to 0 to model electricity markets. Each unit of inflow can be associated with cost if $C^{I}_{v, f, t, \omega} > 0$ or income if $C^{I}_{v, f, t, \omega} < 0$. The parameter $C^{I}_{v, f, t, \omega}$ is included in the cost parameter $c_a$ on all incoming arcs $a$ to demand site $v$. A demand site can also be used to model excess heat by providing an additional heat demand site without any demand (but maybe some penalty costs).

\subsubsection{Storage units $\mathcal{S}$}
Storage units are given in the set $\mathcal{S} \subset \mathcal{V}$ and can store energy from one period to the next. A storage $v \in \mathcal{S}$ is defined for a specific energy type. The conversion factor for this energy type is $\phi_{v,f,f} = (1 - \text{loss})$ to model losses between periods, i.e, the efficiency of the storage unit. For all other energy types, it is set to zero. The capacity of the storage unit is limiting the maximum outflow on arcs connecting the storage with itself in the next period, i.e., $u_a = \text{capacity}(t,\omega)$ with $a = (v,v,f,t,t+1,\omega), \forall t \in \mathcal{T}, \omega \in \Omega$. The capacity can be time- and scenario-dependent. The maximum flow through the storage per period is limiting the total in- and outflow with the lower limits zero, i.e., $[\underline{I}_{v,f,t,\omega}=0,\overline{I}_{v,f,t,\omega}=\text{max-flow}]$ and $[\underline{O}_{v,f,t,\omega}=0,\overline{O}_{v,f,t,\omega}=\text{max-flow}]$. The limitations for all other energy types are zero.

\subsubsection{Interconnections $\mathcal{I}$}
Interconnections are defined by the set $\mathcal{I} \subset \mathcal{V}$ and can be used to model flow restrictions from parts of the network to other parts of the network, e.g., when
% having several distributions grids connected to a transmission grid that has limitations
two districts are connected. An interconnection is given for a defined energy type and for this energy type there are limitations on the in- and outflow given by $[\underline{I}_{v,f,t,\omega},\overline{I}_{v,f,t,\omega}]$ and $[\underline{O}_{v,f,t,\omega},\overline{O}_{v,f,t,\omega}]$, respectively, which can be time- and scenario-dependent. The limitations for all other energy types are zero. Losses can be modelled similar to storages with $\phi_{v,f,f} = (1 - \text{loss})$ for the given energy type.

\subsubsection{Production units $\mathcal{U}$}
Production units are defined by the set $\mathcal{U} \subset \mathcal{V}$. Each production $v \in \mathcal{U}$ can transform input energy types $f_1 \in \mathcal{F}$ to output energy types $f_2 \in \mathcal{F}$. The capacity of the production is unit is defined by the limitations  on the input $[\underline{I}_{v,f1,t,\omega},\overline{I}_{v,f1,t,\omega}]$ and output flows $[\underline{O}_{v,f2,t,\omega},\overline{O}_{v,f2,t,\omega}]$ of each energy type. Non-valid energy types are excluded by setting input and/or output restrictions to zero, respectively. The conversion factor $\phi_{v,f_1,f_2}$ is defined as the relationship between input and output energy type. The increase and decrease in production from one period to the next is limited by the up and down ramping limits $R^U_{v,f}$ and $R^D_{v,f}$, respectively. Those can be defined per energy type $f$. A subset of the units $\mathcal{U}^{*} \subseteq \mathcal{U}$  might relate to here-and-now decisions, i.e., those units need to have the same production in all scenarios $\Omega$ for the given non-anticipativity period $\mathcal{T}^*$.

Some production units require the modelling of their online status (on/off) to impose further restrictions, e.g. in case there exists a minimum production amount or dependencies with other units. The decisions regarding the status of the unit are called commitment decisions and the set of units with commitment decisions is denoted by $\mathcal{U}^{WC}\subseteq \mathcal{U}$. The start-up costs of a unit with commitment are denoted by $C^S_v$. For those units, we can also define minimum up and down times $T^{UT}_v$ and $T^{DT}_v$, respectively. Furthermore, interdependence
between units can be modelled. The set $\mathcal{U}^{DEP}_v$ contains all units that need to be online, when unit $v$ is online.  In contrast, the set $\mathcal{U}^{EXC}_v$ contains all units that are excluded from production, if unit is $v$ producing. This can be used to model one unit with two different operational modes as two units excluding each other.

% \begin{table}
% \caption{Examples}
% \centering
%     \begin{tabular}{p{0.15\columnwidth}p{0.4\columnwidth}p{0.3\columnwidth}}\toprule
%     \textbf{Set} & \textbf{Example} & \textbf{Parameter values}\\\midrule
%     Energy types $\mathcal{F}$ & Natural gas, electricity, heat & $\mathcal{F} = \lbrace NG, EL, HEAT \rbrace$\\\midrule
%     Units $\mathcal{U}$ &  Natural gas-fire CHP unit $u_1$ with max. 20 MWh natural gas as input and 9 MWh heat and 8 MWh power as output  & 
%         $\forall t \in \mathcal{T}, \omega \in \Omega\text{:}$
%         $\overline{I}_{u_1,NG,t,\omega} = 20,$
%         $\overline{O}_{u_1,HEAT,t,\omega} = 9,$
%         $\overline{O}_{u_1,EL,t,\omega} = 8,$ 
%         $\phi_{u_1, NG, EL} = 0.4,$
%         $\phi_{u_1, NG, HEAT} = 0.45$\\
%     Units $\mathcal{U}$ &  Solar plant $u_2$ with max. capacity 1 MWh takes as input solar heat production scenarios and outputs heat  & 
%         $\forall t \in \mathcal{T}, \omega \in \Omega\text{:}$
%         $\overline{I}_{u_2,SolarHeat,t,\omega} = 1,$
%         $\overline{O}_{u_2,HEAT,t,\omega} = 1,$
%         $\phi_{u_1, SolarHeat, Heat} = 1.0$
%         \\
%     \end{tabular}
    
%     \label{tab:my_label}
% \end{table}

\subsection{Network creation}

Based on the structures defined above, the network can be created by translating all components with their parameters to vertices and arcs. Additional information encoded in the network, apart from the attributes mentioned in the previous section, are the limits on particular connections between components based on energy type, time and scenario information. Those are given by bounds on the arcs $[l_a, u_a]$. This means, if two vertices are not connected then the flow bounds are set to zero. The same holds if a particular energy type flow is not possible.  

\subsubsection{Vertices $\mathcal{V}$}
All physical assets in the network, as presented in Section \ref{sec:components}, form the set of vertices in the network, i.e., $\mathcal{V}$ is defined as
$$\mathcal{V} = \mathcal{U} \cup \mathcal{E} \cup \mathcal{D} \cup \mathcal{I} \cup \mathcal{S}.$$
Note that several artificial vertices can be used to account for special conditions. For example, there are artificial energy sources and demand sites for each storage to model initial  levels at the beginning of the planning horizon and target storage levels at the end of the planning horizon. 

\subsubsection{Arcs $\mathcal{A}$}

 The following arcs are created within a period $t \in \mathcal{T}$, $t'=t$ and scenario $\omega \in \Omega$, i.e., flow in the same period and scenario:
\begin{itemize}
\setlength\itemsep{-0.3em}
    \item From energy source $e \in \mathcal{E}$ to unit $u \in \mathcal{U}$, if those two are interconnected and unit $u$ has inflow of the energy type $f$ from source $e$. The upper flow limit is the maximum inflow of energy type $f$ to unit $u$ and the cost are the cost per unit from the energy source.
    \item From energy source $e \in \mathcal{E}$ to storage unit $s \in \mathcal{S}$, demand site $d \in \mathcal{D}$ or interconnection $i \in \mathcal{I}$, if those two are interconnected and end component has inflow of the energy type $f$ of source $e$. The upper limit is the maximum outflow of energy source $e$ and the cost are the cost per unit from the energy source.
     One example is solar heat or waste heat that can be directly used for heating. Then the outflow is limited by the available heat. Another example is the modelling of soft constraints, e.g., imbalances on the electricity markets or missing heat production. In that case, the flow is unlimited and the costs represent penalty costs.
    \item From unit $u \in \mathcal{U}$ to storage unit $s \in \mathcal{S}$, if unit $u$ produces the energy type $f$ stored in $s \in \mathcal{S}$ and the two are connected. The upper limit is the maximum outflow of energy type $f$ to unit $u$.  The cost are the production cost per unit of outflow.
    \item From unit $u \in \mathcal{U}$ to demand site $d \in \mathcal{D}$, if unit $u$ produces the energy type consumed at $d$ and the two are connected. The upper limit is the maximum outflow of energy type $f$ to unit $u$. he cost are the production cost per unit of outflow minus potential income (or cost) for selling the energy type (e.g. electricity on the day-ahead market).
    \item From unit $u_1 \in \mathcal{U}$ to unit $u_2 \in \mathcal{U}$, if $u_1$ produces an energy type $f$ that is input to unit $u_2$ and the two are connected.  The upper limit is the maximum outflow of energy type $f$ to unit $u_1$. The cost are the production cost per unit of outflow.
    \item From unit $u \in \mathcal{U}$ to interconnection $i \in \mathcal{I}$, if unit $u$ produces the energy type of interconnection $i$ and the two are interconnected.  Upper limit is the maximum outflow of energy type $f$ to unit $u$. The cost are the production cost per unit of outflow.
    \item From storage unit $s \in \mathcal{S}$ to demand site $d \in \mathcal{D}$, if storage unit $s$ stores the energy type $f$ of demand site $d$ and the two are connected. The upper limit is the maximum flow per period of storage $s$. The cost are the cost (or income) for sending one unit to the demand site.
    \item From storage unit $s \in \mathcal{S}$ to interconnection $i \in \mathcal{I}$ or vice versa, if storage unit $s$ stores the energy type of interconnection $i$ and the two are connected.  The upper limit is the maximum flow per period of storage $s$. 
    \item From interconnection $i \in \mathcal{I}$ to demand site $d \in \mathcal{D}$, if the interconnection transports the energy type of the demand site and the two are connected.  The upper limit is the maximum flow per period of $i$. The cost are the cost (or income) for sending one unit to the demand site.
    \item From interconnection $i_1 \in \mathcal{I}$ to  interconnection $i_2 \in \mathcal{I}$, if the interconnections  transport the same energy type and the two are connected.  The upper limit is the maximum flow. 
    
\end{itemize}  

Furthermore, the modelling of storage units requires the following arcs from period $t \in \mathcal{T}$ to period $t+1 \in \mathcal{T}$ in the same scenario $\omega \in \Omega$:
\begin{itemize}
    \item From storage unit $s \in \mathcal{S}$ to the same storage unit $s$ in the next period. The upper flow limit is the storage capacity.
\end{itemize}

Special arcs are created for initial and end storage levels at the storage units $s \in \mathcal{S}$. Initially stored energy can enter the system via an artificial energy source $e^* \in \mathcal{E}$ and target storage level can leave the system through artificial demand site $d^* \in \mathcal{D}$. The following arcs are added:
\begin{itemize}
\setlength\itemsep{-0.3em}
    \item In period $t=1$ (first period) from energy source $e^*$ to storage unit $s$ with a lower and upper limit equal to the initial storage level.
    \item In period $t=|\mathcal{T}|$ from storage unit $s$ to demand site $d^*$ with the target storage level as lower bound and storage capacity or target storage level as upper bound.
\end{itemize}

The example in Figure \ref{fig:vertex} could be a setup where vertex $v$ is an electric boiler that can draw electricity ($f$) from either the market ($e_2$) or a wind farm ($u$) that is dependent on wind ($e_1$). The electric boiler produces heat ($f'$) for demand sites $d_1$ and $d_2$ and the flow to demand site $d_2$ is going through a thermal storage ($s$). $\phi_{v,f,f'}$ is the transformation factor from electricity ($f$) to heat ($f'$) for the electric boiler ($v$).

\section{Network flow formulation}\label{sec:model}

\begin{table}
\footnotesize
\centering
	\caption{Sets, parameters and variables}
	\begin{tabular}{lp{0.8\columnwidth}}
		\toprule
		\multicolumn{2}{l}{Variables}\\\midrule
		$x_{a} \in \mathbb{R}^+$ & Flow on arc $a \in \mathcal{A}$\\
% 		$y_{j,t}$ & Binary variable, 1 if unit $j \in \mathcal{V}^\text{C1}$ is producing in period $t$, 0 otherwise\\
% 		$y^{S}_{j,t}$ & Binary variable, 1 if unit $j \in \mathcal{V}^\text{C1}$ is started in period $t$, 0 otherwise\\
% 		$y^{E}_{j,t}$ & Binary variable, 1 if unit $j \in \mathcal{V}^\text{C1}$ is stopped in period $t$, 0 otherwise\\
		$z_{u,t,\omega} \in \{0, 1\}$ & Binary variable, 1 if unit $u \in \mathcal{U}^\text{WC}$ is producing in period $t$ and scenario $\omega \in \Omega$, 0 otherwise\\
		$z^{S}_{u,t,\omega} \in \{0, 1\}$ & Binary variable, 1 if unit $u \in \mathcal{U}^\text{WC}$ is started in period $t$ and scenario $\omega \in \Omega$, 0 otherwise\\
		$z^{E}_{u,t,\omega} \in \{0, 1\}$ & Binary variable, 1 if unit $u \in \mathcal{U}^\text{WC}$ is shut down in period $t$ and scenario $\omega \in \Omega$, 0 otherwise\\\bottomrule

	\end{tabular}
    \label{tab:dv}
\end{table}

{\allowdisplaybreaks
Based on the above defined network, we can create the optimization model. The main decisions of the model are represented by the flow $x_a$ on the arcs $a \in \mathcal{A}$. Further decisions are related to the commitment status of units $u \in \mathcal{U}^{WC}$, where binary variables $z_{u,t,\omega}, z^S_{u,t,\omega}$ and $z^E_{u,t,\omega}$ model whether unit $u$ is online, started up or shut down in period $t$ and scenario $\omega$, respectively. See Table \ref{tab:dv} for an overview of all decision variables including the ranges.

The objective function \eqref{eq:obj} minimizes flow costs through the network plus unit start-up costs.
\begin{align}
	&\min \qquad \sum_{a \in \mathcal{A}} c_{a} x_{a} + \sum_{u \in \mathcal{U}^{WC}}\sum_{t\in T} \sum_{\omega \in \Omega} \pi_\omega C^S_{u} z^{S}_{u,t,\omega}\label{eq:obj}
\end{align}
Constraints \eqref{eq:arc_bounds} limits the flow on the arcs depending on the given bounds.
\begin{align}
	& l_{a} \le x_{a} \le u_{a} && \forall a \in \mathcal{A}\label{eq:arc_bounds}\end{align}
The transformation from one energy type $f$ to another energy $f'$ at a vertex $k$ is handled in constraints \eqref{eq:transformation} using the transformation factor $\phi_{u,f,f'}$. These constraints hold not for energy sources $\mathcal{E}$ and demand sites $\mathcal{D}$, since they do not have incoming or outgoing flow, respectively.
\begin{align}
	 &\sum_{a \in \mathcal{A}^{OUT}_{v,f',t,\omega}} \phi_{v,f,f'} x_{a} &- \sum_{a \in \mathcal{A}^{IN}_{v,f,t,\omega}} x_{a} = 0 \label{eq:transformation}&& \forall v \in \mathcal{V} \backslash (\mathcal{E} \cup \mathcal{D}), t \in \mathcal{T}, \omega \in \Omega, f \in \mathcal{F}, f' \in \mathcal{F}
\end{align}
The total outflow and inflow at each vertex is limited by lower and upper bounds in constraints \eqref{eq:outflow_V} to \eqref{eq:inflow_V}. An exception is made for units with commitment decisions $\mathcal{U}^{WC}$, since those are handled explicitly in constraints \eqref{eq:outflow_VC2} and \eqref{eq:inflow_VC2}.
\begin{align}
	&\underline{O}_{v,f,t,\omega}\le  \sum_{a \in \mathcal{A}^{OUT}_{v,f,t,\omega}}  x_{a} \le \overline{O}_{v,f,t,\omega} && \forall v \in \mathcal{V}^{} \backslash \mathcal{U}^{WC}, t \in \mathcal{T}, \omega \in \Omega, f \in \mathcal{F}\label{eq:outflow_V}\\
	&\underline{I}_{v,f,t,\omega} \le \sum_{a \in \mathcal{A}^{IN}_{v,f,t,\omega}} x_{a} \le \overline{I}_{v,f,t,\omega} && \forall v \in \mathcal{V}^{}  \backslash \mathcal{U}^{WC}, t \in \mathcal{T}, \omega \in \Omega, f \in \mathcal{F}\label{eq:inflow_V}
% 		&y_{k,t}\underline{O}_{k,f,t,\omega} \le\hspace{-15pt} \sum_{(k,j,f,t,t',\omega) \in \mathcal{A}}\hspace{-15pt} x_{(k,j,f,t,t',\omega)} \le y_{k,t}\overline{O}_{k,f,t,\omega} && \forall k \in \mathcal{V}^{C1}, t \in \mathcal{T}, \omega \in \Omega, f \in \mathcal{F}\label{eq:outflow_VC1}\\
% 	&y_{k,t}\underline{I}_{k,f,t,\omega} \le \hspace{-15pt}\sum_{(i,k,f,t,t',\omega) \in \mathcal{A}}\hspace{-15pt} x_{(i,k,f,t,t',\omega)} \le y_{k,t}\overline{I}_{k,f,t,\omega} && \forall k \in \mathcal{V}^{C1}, t \in \mathcal{T}, \omega \in \Omega, f \in \mathcal{F}\label{eq:inflow_VC1}\\
		\end{align}
		The ramping constraints for production units, i.e., the allowed difference in production between periods,  for the different units are given in constraints \eqref{eq:rampup_c0} to \eqref{eq:rampdown_c0} except for units with commitment decisions $\mathcal{U}^{WC}$ (see constraints \eqref{eq:rampup_c2}-\eqref{eq:rampdown_c2}).
		\begin{align}
&\sum_{a \in \mathcal{A}^{OUT}_{u,f,t,\omega}} x_{a} - \sum_{a \in \mathcal{A}^{OUT}_{u,f,t-1,\omega}} x_{a}
    \le R^U_{u,f} &&  \forall u \in \mathcal{U} \backslash \mathcal{U}^{WC}, f \in \mathcal{F}, t \in \mathcal{T}, \omega \in \Omega\label{eq:rampup_c0}\\
        &-\sum_{a \in \mathcal{A}^{OUT}_{u,f,t,\omega}} x_{a} + \sum_{a \in \mathcal{A}^{OUT}_{u,f,t-1,\omega}} x_{a}
    \le R^D_{u,f} && \forall u \in \mathcal{U} \backslash \mathcal{U}^{WC}, f \in \mathcal{F}, t \in \mathcal{T}, \omega \in \Omega\label{eq:rampdown_c0}
    \end{align}
    Please note that the $\sum_{a \in \mathcal{A}^{OUT}_{u,f,t-1,\omega}} x_{a}$ in period $t=1$ refers to the initial production level given as input parameter for each unit $\mathcal{U}$.

\subsection{Units with commitment decisions}
In case of units with commitment decisions, we need to impose additional constraints. The commitment status of the unit $z_{u,t,\omega} \in \{0,1\}$ (1=on, 0=off) impacts the production of the unit. The status variable $z_{u,t,\omega}$ is updated using binary variables for starting $z^S_{u,t\omega}\in \{0,1\}$ and stopping $z^E_{u,t\omega}\in \{0,1\}$ the unit. 

Constraints \eqref{eq:commitment1_c2} to \eqref{eq:depend_c2} model commitment related restrictions for units in set $\mathcal{U}^{WC}$. Constraints \eqref{eq:commitment1_c2} and \eqref{eq:commitment2_c2} ensure that the status of the unit is set correctly based on starting and stopping the unit while excluding simultaneous starts and stops. $z_{u,t-1,\omega}$ in period $t=1$ refers to the initial status given as parameter $B_{u,\omega}$ for each unit $\mathcal{U}^{WC}$.  Minimum up- and down-times are modelled in constraints \eqref{eq:commitment1_init} to \eqref{eq:mindowntime_c2}. Constraints \eqref{eq:commitment1_init} sets the status based on the initial status $B_{u}$ and the required remaining periods $T^B_u$ in this status due to minimum up- or down-time. Constraints \eqref{eq:minuptime_c2} to \eqref{eq:mindowntime_c2}  ensure the minimum up- and down-times for the remaining periods, respectively.  Constraints \eqref{eq:exclude_c2} to \eqref{eq:exclude2_c2} exclude the simultaneous production of two units that should not run at the same time.  The opposite case, where simultaneous production is required, is modelled in Constraints \eqref{eq:depend_c2}. 
\begin{align}
% 	&y^{S}_{k,t} - y^{E}_{k,t} = y_{k,t} - y_{k,t-1} && \forall k \in \mathcal{V}^{C1}, t \in \mathcal{T}\label{eq:commitment_c1}\\
% 		&y^{S}_{k,t} + y^{E}_{k,t} \le 1 && \forall k \in \mathcal{V}^{C1}, t \in \mathcal{T}\label{eq:commitment2_c1}\\
		&z^{S}_{u,t,\omega} - z^{E}_{u,t,\omega} = z_{u,t,\omega} - z_{u,t-1,\omega} && \forall u \in \mathcal{U}^{WC}, t \in \mathcal{T}, \omega \in \Omega\label{eq:commitment1_c2}\\
			&z^{S}_{u,t,\omega}+ z^{E}_{u,t,\omega} \le 1 && \forall u \in \mathcal{U}^{WC}, t \in \mathcal{T}, \omega \in \Omega\label{eq:commitment2_c2}\\
				&z_{u,t,\omega} = B_{u} &&\forall u \in \mathcal{U}^{WC}, t \in \lbrace 0,..., T^{B}_u\rbrace, \omega \in \Omega\label{eq:commitment1_init}\\
				% &\sum_{t'=max\lbrace 1, t-T^{UT}_k\rbrace}^{t} y^{S}_{k,t} \le y_{k,t} && \forall k \in \mathcal{V}^{WC}, t \in \mathcal{T} \label{eq:minuptime_c1}\\
			&\sum_{t'=max\lbrace 1, t-T^{UT}_{u}\rbrace}^{t} z^{S}_{u,t',\omega} \le z_{u,t,\omega} && \forall u \in \mathcal{U}^{WC}, t \in \lbrace {T}^{B}_u, \ldots, |\mathcal{T}| \rbrace, \omega \in \Omega \label{eq:minuptime_c2}\\
% 					&\sum_{t'=max\lbrace 1, t-T^{DT}_k\rbrace}^{t} y^{E}_{k,t} \le 1-y_{k,t} && \forall k \in \mathcal{V}^{C1}, t \in \mathcal{T} \label{eq:mindowntime_c1}\\
			&\sum_{t'=max\lbrace 1, t-T^{DT}_u\rbrace}^{t} z^{E}_{u,t',\omega} \le 1-z_{u,t,\omega} && \forall u \in \mathcal{U}^{WC}, t \in \lbrace {T}^{B}_u, \ldots, |\mathcal{T}| \rbrace, \omega \in \Omega \label{eq:mindowntime_c2}\\
% 		& y_{i,t} + y_{j,t} \le 1 && \forall i \in \mathcal{V}^{C1}, j \in \mathcal{V}^{EXC}_i, t \in \mathcal{T}\label{eq:exclude_c1}\\
% 		& y^{S}_{i,t} + y^{E}_{j,t} \le 1 && \forall i \in \mathcal{V}^{C1}, j \in \mathcal{V}^{EXC}_i, t \in \mathcal{T}\label{eq:exclude2_c1}\\
	& z_{u,t,\omega} + z_{u',t,\omega} \le 1 && \negthickspace\forall u \in \mathcal{U}^{WC}, u' \in \mathcal{U}^{EXC}_u, t \in \mathcal{T},\omega \in \Omega\label{eq:exclude_c2}\\
	& z^{S}_{u,t,\omega} + z^{E}_{u',t,\omega} \le 1 && \negthickspace \forall u \in \mathcal{U}^{WC}, u' \in \mathcal{U}^{EXC}_u, t \in \mathcal{T},\omega \in \Omega\label{eq:exclude2_c2}\\
% 			& y_{i,t} = y_{j,t} && \forall i \in \mathcal{V}^{WC}, j \in \mathcal{V}^{DEP}_i, t \in \mathcal{T}\label{eq:depend_c1}\\
		& z_{u,t,\omega} = z_{u',t,\omega} &&  \negthickspace\forall u \in \mathcal{U}^{WC}, u' \in \mathcal{U}^{DEP}_u, t \in \mathcal{T}, \omega \in \Omega\label{eq:depend_c2}
\end{align}

The inflow and outflow restrictions \eqref{eq:outflow_VC2} and \eqref{eq:inflow_VC2} as well as the ramping of the production  \eqref{eq:rampup_c2} and \eqref{eq:rampdown_c2} of the units with commitment decisions are modelled dependent on the status of the unit. 
\begin{align}
    	&z_{u,t,\omega}\underline{O}_{u,f,t,\omega} \le \sum_{a \in \mathcal{A}^{OUT}_{u,f,t,\omega}}x_{a} \le z_{u,t,\omega}\overline{O}_{u,f,t,\omega} && \forall u \in \mathcal{U}^{WC}, t \in \mathcal{T}, \omega \in \Omega, f \in \mathcal{F}\label{eq:outflow_VC2}\\
	&z_{u,t,\omega}\underline{I}_{u,f,t,\omega} \le \sum_{a \in \mathcal{A}^{IN}_{u,f,t,\omega}} x_{a} \le z_{u,t,\omega}\overline{I}_{u,f,t,\omega} && \forall u \in \mathcal{U}^{WC}, t \in \mathcal{T}, \omega \in \Omega, f \in \mathcal{F}\label{eq:inflow_VC2}
		\end{align}
	\begin{align}
&    \sum_{a \in \mathcal{A}^{OUT}_{u,f,t,\omega}} x_{a} - \sum_{a \in \mathcal{A}^{OUT}_{u,f,t-1,\omega}} x_{a}
    \le R^U_{u,f} z_{u,t-1,\omega} + \underline{O}_{u,f,t,\omega}\nonumber z^{Start}_{u,t,\omega} \\&\hspace{0.5\columnwidth} \forall u \in \mathcal{U}^{WC}, f \in \mathcal{F}, t \in \mathcal{T}, \omega \in \Omega\label{eq:rampup_c2}\\
% &    -\sum_{(i,k,f,t,t',\omega) \in \mathcal{A}} x_{(i,k,f,t,t',\omega)} + \sum_{(i,k,f,t-1,t'',\omega) \in \mathcal{A}} x_{(i,k,f,t-1,t'',\omega)}
%     \le R^D_{u,f} y_{u,t} + \underline{O}_{u,f,t,\omega} y^{Stop}_{u,t} && \forall k \in \mathcal{V}^{C1}, ff \in \mathcal{F}, t \in \mathcal{T}, \omega \in \Omega\label{eq:rampdown_c1}\\
&    -\sum_{a \in \mathcal{A}^{OUT}_{u,f,t,\omega}} x_{a} + \sum_{a \in \mathcal{A}^{OUT}_{u,f,t-1,\omega}} x_{a}
    \le R^D_{u,f} z_{u,t,\omega} + \underline{O}_{u,f,t,\omega} z^{Stop}_{u,t,\omega}\nonumber \\&\hspace{0.5\columnwidth}\forall u \in \mathcal{U}^{WC}, f \in \mathcal{F}, t \in \mathcal{T}, \omega \in \Omega\label{eq:rampdown_c2}
\end{align}

\subsection{Non-anticipativity constraints}
Depending on the application case of the model, different non-anticipativity constraints need to be added. If a deterministic version of the model is used, those constraints can be omitted. Furthermore, we can distinguish between non-anticipativity on the commitment decisions of the units and non-anticipativity with respect to bidding curves. The former assumes that the units can operate with the day-ahead market prices without any bidding, while the latter creates bidding curves that can be submitted to the day-ahead electricity market.

\subsubsection{Operational planning without bidding}
For all units with here-and-now decisions $\mathcal{U}^{*}$ in the periods considered as first-stage $\mathcal{T}^*$, we need to include non-anticipativity constraints to ensure the decision structure of a two-stage stochastic program. Non-anticipativity means that the production and status of those units needs to be equal across scenarios. Such non-anticipativity might be necessary, e.g., due to electricity market participation, where the determining of power production amounts has to be made before the market is cleared. The constraints for commitment status and production are given in constraints \eqref{eq:nonanticommit} and \eqref{eq:nonantiflow}, respectively. Constraints \eqref{eq:nonanticommit} ensures that the commitment status is the same for all scenarios for all units in $\mathcal{U}^* \cap \mathcal{U}^{WC}$. The flow non-anticipativity is modeled such that all arcs need to have the expected flow \eqref{eq:nonantiflow}. The set $\mathcal{A}^{*}(a)$ contains all arcs $a'$ that need to have same flow as arc $a$, i.e., the arcs with same start and end vertex, period and energy type, which only differ in scenario, including $a$ itself. The non-anticipativity on flow holds only for arcs with units $v(a)$ in $\mathcal{U}^*$ as start vertex and period $t(a) \in \mathcal{T}^*$ as start period.  ${\omega(a')}$ denotes the scenario of arc $a'$.
	\begin{align}
		&z_{u,t,\omega} = z_{u,t,\omega'} && \forall u \in \mathcal{U}^* \cap \mathcal{U}^{WC}, t \in \mathcal{T}^*, \omega, \omega' \in \Omega\label{eq:nonanticommit}\\
		&x_{a} = \sum_{a' \in \mathcal{A}^{*}(a)}\pi_{\omega(a')} x_{a'} && \forall \lbrace a \in  \mathcal{A}\ |\ v(a) \in \mathcal{U}^* \land t(a) \in \mathcal{T}^* \rbrace\label{eq:nonantiflow}
	\end{align}
}

In this setting (in contrast to the approach described in Section \ref{sec:bidding_curves}) the DH system is assumed to participate in electricity trading via price-independent bids. That means that electricity trades, if first-stage decisions, are realized no matter the market price.

\subsubsection{Operational planning including bidding curves}\label{sec:bidding_curves}
The extension to bidding curves is based on the work in \cite{hvidesande} that use the method of \cite{pandvzic2013offering} for creating bidding curves based on electricity price scenarios. The method creates monotonously increasing/decreasing bidding curves determining a bidding amount for each price scenario. We refer to \cite{hvidesande} and \cite{pandvzic2013offering} for further details.

\begin{figure}
    \centering
    \includegraphics[width=0.6\columnwidth]{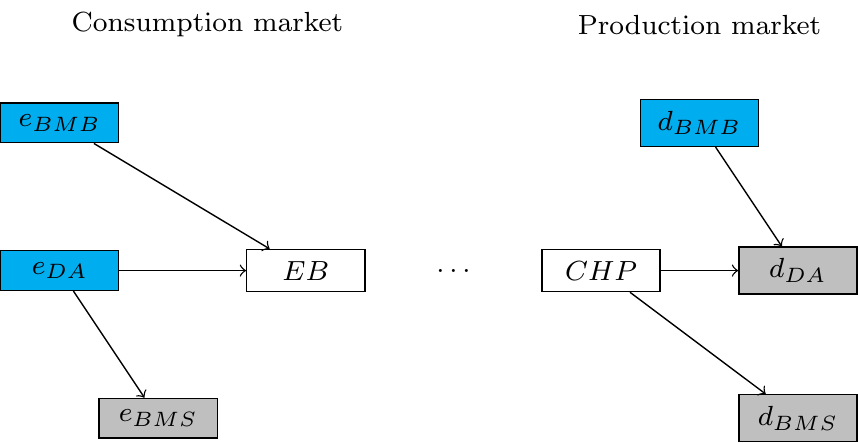}
    \caption{Set of vertices needed to represent markets with exemplary connections to electric boiler (EB) and CHP unit (CHP). Each market consists of three vertices: day-ahead market (DA), negative imbalance (buying needed) (BMB) and positive imbalance (selling needed) (BMS) (blue = energy sources, gray = demand sites, white = units).}
    \label{fig:markets}
\end{figure}

In this extension, we have a set of markets for selling $m \in \mathcal{M}^S$ and buying $m \in \mathcal{M}^B$ energy. A market $m$ contains a set of three vertices $\mathcal{V}_m$ representing a spot (day-ahead) market as well as imbalances (upward and downward). These vertices are either energy sources, if they offer energy, or demand sites, if they receive energy. The price on the spot market $m \in \mathcal{M}^{S} \cup \mathcal{V}^{S}$ in period $t$ and scenario $\omega$ is denoted by $p_{m,t,\omega}$. The setup of vertices is visualized in Figure \ref{fig:markets}.

Constraints \eqref{eq:bidcurvesell} and \eqref{eq:bidcurvesell2} determine the bidding amount to the selling market based on the scenarios and allow a monotonously increasing bidding curve based on the market scenario price $p_{m,t,\omega}$. For an equal price an equal selling amount is guaranteed \eqref{eq:bidcurvesell2} ensuring non-anticipativity. Constraints \eqref{eq:bidcurvebuy} and \eqref{eq:bidcurvebuy2} create monotonously decreasing bidding curves for the buying market.
{\allowdisplaybreaks
\begin{align}
    \sum_{a \in \mathcal{A}^{IN}_{v,f,t\omega}} x_a - \sum_{a \in \mathcal{A}^{OUT}_{v,f,t\omega}} x_a &= \sum_{a \in \mathcal{A}^{IN}_{v,f,t\omega'}} x_a - \sum_{a \in \mathcal{A}^{OUT}_{v,f,t\omega'}} x_a  \nonumber\\&\forall m \in \mathcal{M}^S, v \in \mathcal{V}_m, t \in \mathcal{T}, (\omega, \omega') \in \Omega \times \Omega, \text{if } p_{m,t,\omega} = p_{m,t\omega'}\label{eq:bidcurvesell}\\
    \sum_{a \in \mathcal{A}^{IN}_{v,f,t\omega}} x_a - \sum_{a \in \mathcal{A}^{OUT}_{v,f,t\omega}} x_a &\le \sum_{a \in \mathcal{A}^{IN}_{v,f,t\omega'}} x_a - \sum_{a \in \mathcal{A}^{OUT}_{v,f,t\omega'}} x_a \nonumber\\& \forall m \in \mathcal{M}^S, v \in \mathcal{V}_m, t \in \mathcal{T}, (\omega, \omega') \in \Omega \times \Omega, \text{if } p_{m,t,\omega} \le p_{m,t\omega'}\label{eq:bidcurvesell2}\\
    \sum_{a \in \mathcal{A}^{OUT}_{v,f,t\omega}} x_a - \sum_{a \in \mathcal{A}^{IN}_{v,f,t\omega}} x_a &= \sum_{a \in \mathcal{A}^{OUT}_{v,f,t\omega'}} x_a - \sum_{a \in \mathcal{A}^{IN}_{v,f,t\omega'}} x_a  \nonumber\\&\forall m \in \mathcal{M}^B, v \in \mathcal{V}_m, t \in \mathcal{T}, (\omega, \omega') \in \Omega \times \Omega, \text{if } p_{m,t,\omega} = p_{m,t\omega'}\label{eq:bidcurvebuy}\\
    \sum_{a \in \mathcal{A}^{OUT}_{v,f,t\omega}} x_a - \sum_{a \in \mathcal{A}^{IN}_{v,f,t\omega}} x_a &\le \sum_{a \in \mathcal{A}^{OUT}_{v,f,t\omega'}} x_a - \sum_{a \in \mathcal{A}^{IN}_{v,f,t\omega'}} x_a \nonumber\\& \forall m \in \mathcal{M}^B, v \in \mathcal{V}_m, t \in \mathcal{T}, (\omega, \omega') \in \Omega \times \Omega, \text{if } p_{m,t,\omega} \ge p_{m,t\omega'}\label{eq:bidcurvebuy2}
\end{align}}

\subsection{Rolling horizon approach}
The model presented above can also be used in a rolling horizon setting where we shift the planning horizon by  $|\mathcal{T}^*|$ in each iteration. To model the rolling horizon correctly, some input parameters need to be updated based on the realization of the uncertainty. The initial status of the units with commitment decisions $\mathcal{U}^{WC}$ and the initial production and storage levels need to be set according to the outcome in period $t=|\mathcal{T}^*|$ in the previous run.

\section{Test cases}\label{sec:cases}

We present results for the optimization of three different DH systems in Denmark. These are located in the cities of Brønderslev, Hillerød and Middelfart.

\subsection{Static data}
\begin{figure}
     \centering
     \begin{subfigure}[t]{0.49\columnwidth}
         \centering
         \includegraphics[width=0.9\columnwidth]{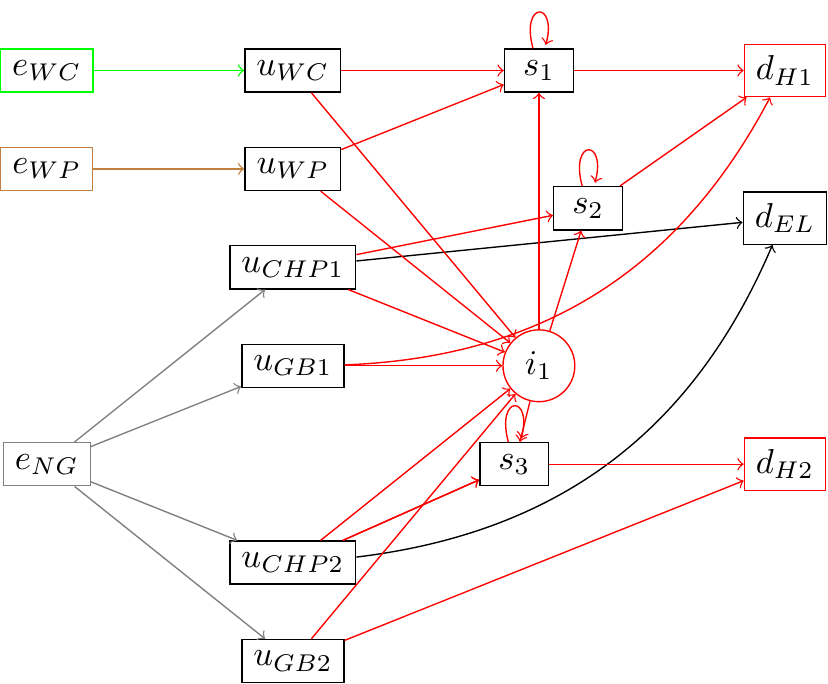}
         \caption{Middelfart}
         \label{fig:middelfart}
     \end{subfigure}
     \hfill
     \begin{subfigure}[t]{0.49\columnwidth}
         \centering
         \includegraphics[width=0.9\columnwidth]{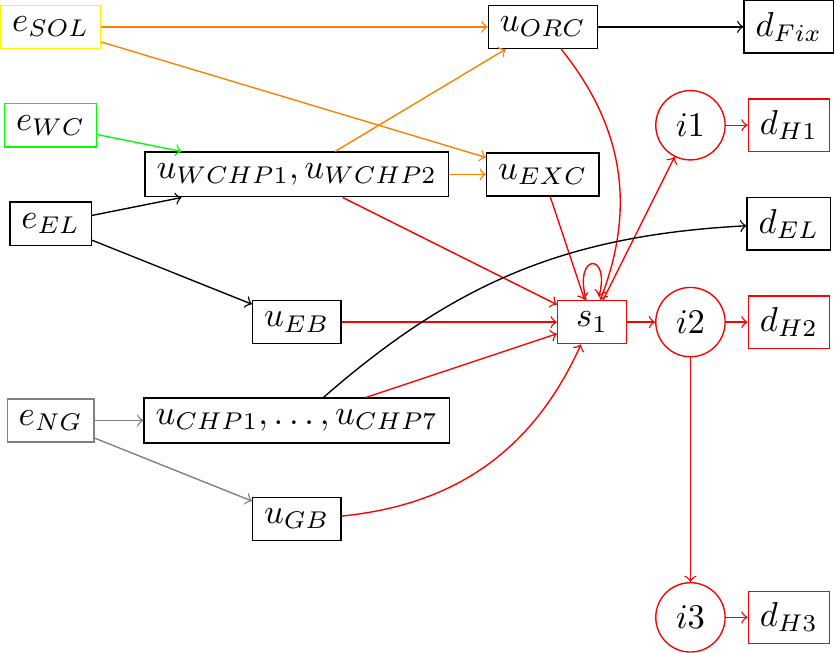}
         \caption{Brønderslev}
         \label{fig:bronderslev}
     \end{subfigure}
     \hfill
     \begin{subfigure}[t]{\columnwidth}
         \centering
         \includegraphics[width=0.5\columnwidth]{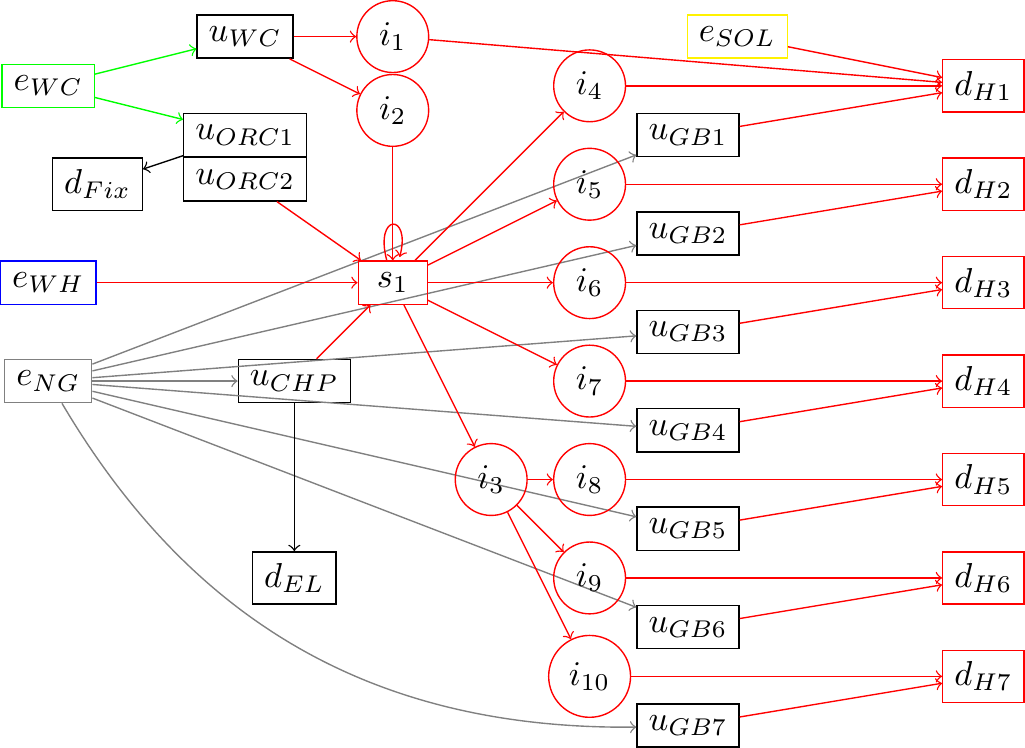}
         \caption{Hillerød}
         \label{fig:hilleroed}
     \end{subfigure}
        \caption{Basic network representation of the three DH systems with energy sources ($e$), units ($u$), storage units ($s$), demand sites ($d$) and interconnections ($i$). For the Brønderslev system, the seven CHP units and two wood chip boiler-heat pump units are aggregated for illustrative purposes. WC=wood chips, WP=wood pellets, CHP=combined head and power, GB=gas boiler, EL=electricity, H=heat, SOL=solar heat, ORC=organic rankine cycle, EXC=heat exchanger, NG=natural gas.}
        \label{fig:dhs}
\end{figure}

The DH systems differ in terms of types and numbers of units as well as the layout of the demand sites. An overview of the basic network representation is given in Figure \ref{fig:dhs}. 

% Table generated by Excel2LaTeX from sheet 'Sheet2'
\begin{table}[t]
  \centering
  \footnotesize
  \caption{Unit parameters ($^+$ = Cost and prices have been multiplied with a factor to anonymize data, $^*$ = First-stage decisions, i.e., $u \in \mathcal{U}^*$)}
  \begin{adjustbox}{width=0.9\textwidth}
    \begin{tabular}{p{0.01\columnwidth}p{0.06\columnwidth}p{0.2\columnwidth}p{0.23\columnwidth}p{0.1\columnwidth}p{0.2\columnwidth}}
          \toprule
          & Unit  & Input restr.  & Output restr.  & \multicolumn{1}{l}{Cost } & Min. up/down  \\
          &       & [MW/period] & [MW/period] & [EUR/
          MWh$_h$] &times [periods], \qquad Start up cost [EUR] \\\midrule
    \multirow{10}[0]{*}{\begin{sideways}Brønderslev\end{sideways}} & \multicolumn{1}{l}{$u_{CHP1}^*$-} & {NG: [7.822,7.822]} & H: [3.848, 3.848],  & {96.98} & \multicolumn{1}{l}{{}} \\
          & \multicolumn{1}{l}{$u_{CHP7}^*$} &       & EL: [3.277, 3.277] &       &  \\
          & $u_{GB}$    & NG: 33.8 & H: 34.8 & 58.52 &  \\
          & $u_{SOL}$   & SOL: External & H: External & 0.00  &  \\
          & $u_{ORC}$   & PH: [10, 20] & H: [8,16], EL: [2,4] & -32.16 & \multicolumn{1}{l}{Up/Down: 3/3} \\
          & $u_{EB}^*$    & EL: 20 & H: 19.8 & 56.88 & {} \\
          & {$u_{WCHP1}^*$} & WC: [4.32, 10.8]  & H: [1.2, 3],  & {88.51} & {{}} \\
          &       & EL: [0.12, 0.3] & PH: [3.2, 8] &       &  \\
          & $u_{WCHP2}^*$ & WC: 10.8 EL: 0.3 & H: 3 PH: 8 & 88.51 & {dep. on WCHP1} \\
          & $u_{EXC}$   & PH: 20 & H: 20 & 0.00  &  \\\midrule
    \multirow{15}[0]{*}{\begin{sideways}Hillerød$^+$\end{sideways}} & $u_{WC}$    & WC: 7.9 & H: 8  & 33.74 & {Start up cost: 30} \\
          & {$u_{CHP}^*$} & {NG: [58, 146]} & H: [26, 65],  & {100} & {Start up cost: 500} \\
          &       &       & EL: [13, 59.1] &       &  \\
          & $u_{GB1}$   & NG: 28.4 & H: 27 & 58.19 &  \\
          & $u_{GB2}$   & NG: 4.6 & H: 4.5 & 58.22 &  \\
          & $u_{GB3}$   & NG: 2 & H: 1.9 & 58.22 &  \\
          & $u_{GB4}$   & NG: 16.842 & H: 16 & 58.22 &  \\
          & $u_{GB5}$   & NG: 23.368 & H: 22.2 & 58.22 &  \\
          & $u_{GB6}$   & NG: 35.79 & H: 34 & 56.69 &  \\
          & $u_{GB7}$   & NG: 20.4 & H: 20 & 56.58 &  \\
          & $u_{SOL}$   & SOL: External & H: External & 0.00  &  \\
          & $u_{ORC1}^*$  & WC: [6.86, 30.83] & H: [5, 25.6], EL: [0, 4] & 43.57 & {Up/Down: 24/24} \\
          & {$u_{ORC2}^*$} & {WC: [6.86, 30.83]} & H: [4.675, 23.936],  & {44.32} & {Up/Down: 24/24,} \\
          &       &       & EL: [0, 3.740] &       & {excludes ORC1} \\
          & $u_{WH}$    & WH: External & H: External & 0.00  &  \\\midrule \multirow{8}[0]{*}{\begin{sideways}Middelfart\end{sideways}} & $u_{WC}$    & WC: [0.775, 4.095] & {H: [0.814; 4.3]} & 24.19 & \multicolumn{1}{l}{Up/Down: 24/24} \\
          & $u_{WP}$    & WP: [0.555, 2.760] & {H: [0.52, 2.5]} & 30.24 & \multicolumn{1}{l}{Up/Down: 12/12} \\
          & {$u_{CHP1}^*$} & {NG: [7.565, 7.565]} & {H: [3.625, 3.625], } & {109.61} & \multicolumn{1}{l}{Start up cost: 72.67} \\
          &       &       & {EL: [2.875, 2.875]} &       & \multicolumn{1}{l}{} \\
          & {$u_{CHP2}^*$} & {NG: [7.675, 7.675]} & {H: [4.22, 4.22], } & {64.13} & \multicolumn{1}{l}{{Start up cost: 73.72}} \\
          &       &       & {EL: [3.3, 3.3]} &       &  \\
          & $u_{GB1}$  & NG: 5.59 & {H: 5.815} & 63.08 &  \\
          & $u_{GB2}$  & NG: 6.33 & {H: 6.52} & 46.67 &  \\\bottomrule
    \end{tabular}%
\end{adjustbox}
  \label{tab:unitparams}%
\end{table}%

% Table generated by Excel2LaTeX from sheet 'Sheet1'
\begin{table}[t]
  \centering
  \footnotesize
  
  \caption{Lower and upper bounds as well as prices of energy sources ($e$) and demand sites ($d$). TS=time series, Scen.=Scenario set ($^*=$ the price is indirectly taken into account since it is included in the heat production cost of the unit)}
    \begin{adjustbox}{width=0.8\textwidth}
    \begin{tabular}{llll}\toprule
          & \multicolumn{1}{l}{LB [MW/period]} & UB [MW/period] & Price [EUR/MWh] \\\midrule
    $e_{WC}$ & 0     & Unlimited & 0$^*$ \\
    $e_{WP}$ & 0     & Unlimited & 0$^*$ \\
    $e_{NG}$ & 0     & Unlimited & 0$^*$ \\
    $e_{EL}$ & 0     & Unlimited & Day-ahead price TS/Scen \\
    $e_{SOL}$ & \multicolumn{1}{l}{Solar production TS/Scen} & Solar production TS/Scen & \multicolumn{1}{l}{0} \\
    $e_{WH}$ & \multicolumn{1}{l}{Waste heat prod. TS/Scen} & Waste heat prod. TS/Scen & \multicolumn{1}{l}{0} \\\midrule
    $d_{H\cdot}$ & \multicolumn{1}{l}{Heat demand TS/Scen} & Heat demand TS/Scen & \multicolumn{1}{l}{0} \\
    $d_{EL}$ & 0     & Unlimited & Day-ahead price TS/Scen \\
    $d_{Fix}$ & 0     & Unlimited & 0$^*$ \\\bottomrule
    \end{tabular}%
    \end{adjustbox}
  \label{tab:tsdata}%
\end{table}%

% Table generated by Excel2LaTeX from sheet 'Sheet1'
% Table generated by Excel2LaTeX from sheet 'Sheet1'
\begin{table}[ht]
\footnotesize
  \centering
  \caption{Interconnection limits for each DH system}
  \begin{adjustbox}{width=0.75\textwidth}
    \begin{tabular}{lrrrrrrrrrrrr}\toprule
          &  \multicolumn{3}{c}{Brønderslev} & \multicolumn{8}{c}{Hillerød} &  \multicolumn{1}{l}{Middelfart} \\\cmidrule(lr){2-4}\cmidrule(lr){5-12} \cmidrule(lr){13-13}
    $i$      & 1     & 2     & 3     & 1     & 2     & 3     & 4     & 5     & 6     & 7     & 8 & 1 \\\midrule
    LB [MW/period]      & 0     & 0     & 0     & 0     & 0     & 0     & 0     & 0     & 0     & 0     & 0 & 0 \\
    UB [MW/period]      & 20    & 32    & 12    & 20    & 5     & 2     & 16.5  & 22.2  & 16    & 2     & 2 & 5\\\bottomrule
    \end{tabular}%
 \end{adjustbox}
  \label{tab:idata}%
\end{table}%

The Middelfart DH system in Figure \ref{fig:middelfart} only contains a few units and consists of two sub-systems connected through an interconnection pipe($i_1$). One sub-system contains a wood chip boiler ($u_{WC}$), wood pellet boiler ($u_{WP}$), a CHP unit  ($u_{CHP1}$) and a gas boiler ($u_{GB1}$), while the other sub-system contains only a gas boiler ($u_{GB2}$) and a CHP unit ($u_{CHP2}$). The entire system contains three storage units ($s_1-s_3$) and two demand sites ($d_{H1}, d_{H2}$). The CHP units sell electricity to the day-ahead market ($d_{el}$).

The Brønderslev DH system in Figure \ref{fig:bronderslev} has a larger number and variety of units. Apart from 7 small-scale CHP units ($u_{CHP1},\ldots,u_{CHP7}$), an electric boiler ($u_{EB}$) and a gas boiler ($u_{GB}$) producing heat used directly in the DH network, the system has two units combining wood chip boilers with heat pumps ($u_{WCHP1}, u_{WCHP2}$) and a solar thermal plant ($u_{SOL}$) to produce process heat (PH). The process heat is then utilized by an organic rankine cycle (ORC) CHP unit ($u_{ORC}$) for combined heat and power production or it can be converted to heat through an heat exchanger ($u_{EXC}$). The heat is delivered to a central storage ($s_1$) and from there, via interconnections $i_{1} - i_{3}$, to three demand sites ($d_{H1} - d_{H3}$) where $d_{H3}$ is reached via $d_{H2}$. The CHP (to the day-ahead market ($d_{EL}$)) and ORC units (at fixed price) sell electricity, while the electric boiler and heat pumps consume electricity from the day-ahead market. 

The Hillerød DH system in Figure \ref{fig:hilleroed} is characterized by eight separate demand sites ($d_{H1} - d_{H8}$) that are reached via interconnections ($i_1 - i_8$), nearly all of them connected to decentralized gas boilers ($u_{GB1} - u_{GB7}$) that can cover local heat production. Demand site $d_{H1}$ is also connected to a solar thermal plant for heat production ($u_{SOL}$). At the central production, waste heat injection ($u_{WH}$), a wood chip boiler ($u_{WC}$), one CHP unit ($u_{CHP}$) and one wood chip-fired ORC unit are able to produce heat. The ORC unit can be operated in two different operational modes represented by units $u_{ORC1}$ and $u_{ORC2}$, production by either excluding production by the other. The electricity production from the CHP unit is sold on the day-ahead market ($d_{EL}$) while the electricity production from the ORC unit is sold at a fixed price ($d_{Fix}$). 

For all systems, the parameters of the units are given in Table \ref{tab:unitparams}, the parameters for energy sources and demand sites are stated in Table \ref{tab:tsdata} and the data for interconnections and storage units are given in Table \ref{tab:idata} and \ref{tab:sdata}, respectively. We use real-world data provided by the respective DH system operator, except all costs and prices for the Hillerød system that were multiplied with a constant factor to anonymize the costs but keep the relation between units and markets.  

The day-ahead market $d_{El}$ in all three systems is represented by day-ahead market and two vertices for imbalances (see description in Section \ref{sec:bidding_curves}). The penalty costs for imbalances are 600 EUR/MWh (higher than all electricity prices in the dataset). Additionally, each system contains an energy source for missing heat (with a penalty cost of 10 EUR/kWh) and a demand site for excess heat (with no penalty).

% Table generated by Excel2LaTeX from sheet 'Sheet1'
\begin{table}
  \centering
  \footnotesize
  \caption{Storage unit parameters for each DH system}
  \begin{adjustbox}{width=0.8\textwidth}
    \begin{tabular}{llllll}\toprule
          & {$s$} & Capacity [MWh] & Initial level [MWh] & Target level [MWh] & Loss \quad[\%/period] \\\midrule
    Brønderslev & 1     & 361.54 & 0.1    & 0.1    & 0.01 \\\midrule
    Hillerød & 1     & 556.22 & 0.1   & 0.1      & 0.01 \\\midrule
     \multirow{3}[0]{*}{Middelfart}  & 1     & 38.048 & 0.1    & 0.1      & 0.01 \\
          & 2     & 47.56 & 0.1    & 0.1     & 0.01\\
          & 3     & 41.136 & 0.1    & 0.1      & 0.01 \\\bottomrule
    \end{tabular}%
    \end{adjustbox}
  \label{tab:sdata}%
\end{table}%

\subsection{Time series data and scenarios}
 For each system, we analyze cases from different seasons and varying planning horizon lengths. Evaluating system behaviour during different periods of the year is important to capture seasonal variations. An overview of test cases is provided in Table \ref{tab:cases}. We divide into stochastic optimization for operational planning and deterministic optimization for longer planning horizons. The purpose of the former is the operational optimization considering the interaction with the day-ahead electricity market including the operational uncertainty. In the latter case, we use deterministic data input to evaluate the operational performance of the system in a particular setting. All datasets have an hourly resolution. Note that the historical data in the three systems is based on different years (2019-2021), i.e., results and costs are not directly comparable due to different weather conditions and market prices.

% Table generated by Excel2LaTeX from sheet 'Sheet2'
% Table generated by Excel2LaTeX from sheet 'Sheet2'
\begin{table}
  \centering
  \footnotesize
  \caption{Test cases}
  \begin{adjustbox}{width=0.8\textwidth}
    \begin{tabular}{llllp{0.15\columnwidth}p{0.05\columnwidth}}\toprule
    Name  & \multicolumn{1}{l}{System} & \multicolumn{1}{l}{Start date} & Dataset length & Planning horizon in one model run &  det./ sto. \\
    \midrule
    B-01-168 & \multirow{3}[2]{*}{Brønderslev} & 18-01-2021 & 2 weeks & 1 week    & sto. \\
    B-04-168 &       & 26-04-2021 & 2 weeks & 1 week    & sto. \\
    B-07-168 &       & 05-07-2021 & 2 weeks & 1 week    & sto. \\
    B-10-6936 &       & 15-10-2020 & to 31-07-2021, 289 days & 289 days     & det. \\
    \midrule
    H-04-168 & \multirow{3}[2]{*}{Hillerød} & 05-04-2021 & 2 weeks & 1 week    & sto. \\
    H-07-168 &       & 05-07-2021 & 2 weeks & 1 week    & sto. \\
    H-10-168 &       & 04-10-2021 &2 weeks &  1 week    & sto. \\
      H-02-5808 &       & 15-02-2021 & to 14-10-2021, 242 days & 242 days     & det. \\\midrule
       M-05-168 & \multirow{3}[2]{*}{Middelfart} & 13-05-2019 & 2 weeks & 1 week    & sto. \\
    M-08-168 &       & 05-08-2019 & 2 weeks & 1 week    & sto. \\
    M-12-168 &       & 21-12-2019 & 2 weeks & 1 week    & sto. \\
     M-03-8784 &       & 03-03-2019 & to 16-12-2019, 289 days & 289 days     & det.\\
    \bottomrule

    \end{tabular}%
    \end{adjustbox}
  \label{tab:cases}%
\end{table}%

For each DH system, we used real historical data to evaluate the performance of the models. For the scenario generation, we follow a very basic procedure based on historical data. For each test case, we use the data from the three previous weeks (before the date stated in Table \ref{tab:cases}) weighted with 0.5, 0.33 and 0.17 giving more weight to the recent observations. We distinguish between uncertainty regarding heat flows (solar production, waste heat injection and heat demand) and electricity prices. We combine the price scenarios with the heat flow scenario resulting in a total of 9 scenarios.  This distinction is necessary for the creating of the bidding curves (see \cite{pandvzic2013offering}), i.e. the model creates bidding curves with three steps. This very basic scenario generation can easily replaced by state-of-the-art probabilistic forecasting methods by just exchanging the input data. Additional data used for the evaluation is stated in the respective results section.

\section{Results}\label{sec:results}
All results are computed on hardware of the DTU Computing Center(DCC)\footnote{https://www.hpc.dtu.dk/} with Intel Xeon Processors 2650v4 2.20GHz using 8 cores and 16 GB RAM. The models are implemented using Python 3.7.10 and PuLP 2.5.1, and solved with Gurobi 9.5.0.
 
\subsection{Value of stochastic solution and out-of-sample performance}
 In the first experiment, we evaluate the operational optimization with and without bidding curves comparing the stochastic programming approach with a simpler deterministic model using the expected value of the uncertain data. 
 
 Tables \ref{tab:vss-nobidding} and \ref{tab:vss-bidding} present the objective values and several solution metrics for each test case using each of the two approaches with and without bidding, respectively. Furthermore, we calculate the value of stochastic solution (VSS) that is defined as the difference between the expectation of the expected value solution (Exp.) and the expected value of the stochastic program (Sto.). The VSS is a standard metric for evaluation of stochastic programming \cite[pp.165]{birge2011introduction}. It evaluates the performance based on the scenarios considered in the model, i.e., 9 scenarios in our case.  The planning horizon of the model is 168h (first week of the respective dataset) and the first-stage decisions relate to either the commitment status of the units (No bidding) or the bidding curves (Bidding). When bidding is considered, the evaluation first determines if the bid to the market would have been successful in a certain scenario by comparing the bidding price with the scenario price. Only if the bid is successful, production is allowed.
 
%  \begin{table}
% \footnotesize
%  \caption{Value of stochastic solution (VSS), i.e., difference between exp. value approach and stochastic programming evaluated over the 9 scenarios included in the models. Exp., Sto. and VSS are given in EUR. }
%     \centering
%      %\begin{adjustbox}{width=\columnwidth}
%     \begin{tabular}{lrrrrrrrrr}\toprule
%      Case & \multicolumn{4}{c}{No bidding} & \multicolumn{4}{c}{Bidding}\\\cmidrule(r){2-5}\cmidrule(r){6-9}
%         & Exp. & Sto. & VSS & VSS [\%] & Exp. & Sto. & VSS & VSS [\%] \\ \midrule
%             B-01-168 & 118831.00 & 111703.76 & 7127.24 & 6.0\% & 144497.58 & 106370.66 & 38126.92 & 26.4\% \\ 
%         B-04-168 & 37993.11 & 36717.83 & 1275.28 & 3.4\% & 66262.43 & 32010.09 & 34252.33 & 51.7\% \\ 
%         B-06-168 & 19861.07 & 15144.78 & 4716.29 & 23.7\% & 33424.11 & 14327.40 & 19096.71 & 57.1\% \\ 
%         H-04-168 & 282758.56 & 279712.80 & 3045.76 & 1.1\% & 288430.29 & 267060.03 & 21370.26 & 7.4\% \\ 
%         H-07-168 & 3158.05 & 1989.38 & 1168.68 & 37.0\% & 29094.95 & -561.41 & 29656.36 & 101.9\% \\ 
%         H-10-168 & -513749.36 & -513749.36 & 0.00 & 0.0\% & -339867.25 & -520133.29 & 180266.04 & 53.0\% \\ 
%         M-05-168 & 16270.59 & 16270.59 & 0.00 & 0.0\% & 16270.50 & 16270.59 & -0.09 & 0.0\% \\ 
%         M-08-168 & 6971.40 & 6971.34 & 0.06 & 0.0\% & 6899.44 & 6862.05 & 37.39 & 0.5\% \\ 
%         M-12-168 & 33097.00 & 33055.46 & 41.54 & 0.1\% & 33602.95 & 32799.75 & 803.20 & 2.4\% \\  \bottomrule
%     \end{tabular}
%     %\end{adjustbox}
%     \label{tab:vss}
% \end{table}

\begin{table}[t]
    \centering
    \footnotesize
    \caption{Performance of expected value approach and stochastic program \textbf{without} consideration of bidding: Objective value (Obj.), net electricity sold on day-ahead market (sales-purchase) (El.sales), Income from electricity market (Income), Heat production (Heat) and cost per MWh. All values are expected values. The VSS is the difference between objective values and VSS[\%] gives the improvement when using the stochastic program. }
    \begin{adjustbox}{width=\textwidth}
    \begin{tabular}{lrrrrr|rrrrr|rr}
     \toprule
        ~ & \multicolumn{5}{c}{Expected value approach}& \multicolumn{5}{c}{Stochastic program} & ~ & ~ \\ \cmidrule(rl){2-6} \cmidrule(rl){7-11}
        Case & Obj. & El.sales & Income & Heat & Cost/ & Obj. & El.sales & Income & Heat & Cost/ & VSS & VSS\\
          & [EUR]  & [MWh$_\text{e}$] & [EUR] & [MWh$_\text{h}$]& MWh$_\text{h}$ & [EUR]  & [MWh$_\text{e}$] & [EUR] & [MWh$_\text{h}$]& MWh$_\text{h}$ & [EUR]& [\%] \\\midrule
        %\rot{Case} & \rot{Objective value [EUR]} & \rot{El. sales[MWh]} &\rot{Market income [EUR]} & \rot{Heat production [MWh]} & \rot{Cost per MWh-heat} & \rot{Objective value [EUR]} & \rot{El. sales[MWh]} & \rot{Market income [EUR]} & \rot{Heat production [MWh]} & \rot{Cost per MWh-heat} & \rot{VSS [EUR]} & \rot{VSS [\%]} \\ 
        B-01-168 & 118831.0 & 1347.6 & 81222.5 & 5118.3 & 23.2 & 111703.8 & 987.2 & 62646.6 & 4957.3 & 22.5 & 7127.2 & 6.0 \\ 
        B-04-168 & 37993.1 & 6.7 & 5185.9 & 3663.1 & 10.4 & 36717.8 & 5.6 & 4900.1 & 3659.0 & 10.0 & 1275.3 & 3.4 \\ 
        B-06-168 & 19861.1 & -67.4 & -3281.9 & 2268.6 & 8.8 & 15144.8 & -71.6 & -3539.3 & 2330.7 & 6.5 & 4716.3 & 23.7 \\ 
        H-04-168 & 282758.6 & 3458.8 & 274766.1 & 8490.2 & 33.3 & 279712.8 & 3380.5 & 270196.6 & 8490.8 & 32.9 & 3045.8 & 1.1 \\ 
        H-07-168 & 3158.1 & 2155.1 & 242157.0 & 2701.5 & 1.2 & 1989.4 & 2181.1 & 246348.6 & 2749.8 & 0.7 & 1168.7 & 37.0 \\ 
        H-10-168 & -513749.4 & 9361.4 & 1544349.4 & 10518.3 & -48.8 & -513749.4 & 9361.4 & 1544349.4 & 10518.3 & -48.8 & 0.0 & 0.0 \\ 
        M-05-168 & 16270.6 & 0.0 & 0.0 & 668.5 & 24.3 & 16270.6 & 0.0 & 0.0 & 668.5 & 24.3 & 0.0 & 0.0 \\ 
        M-08-168 & 6971.4 & 36.3 & 2071.7 & 291.1 & 23.9 & 6971.3 & 36.3 & 2071.7 & 291.1 & 23.9 & 0.1 & 0.0 \\ 
        M-12-168 & 33097.0 & 184.8 & 8981.9 & 1247.6 & 26.5 & 33055.5 & 214.5 & 10105.9 & 1247.6 & 26.5 & 41.5 & 0.1 \\ \bottomrule
    \end{tabular}
        \label{tab:vss-nobidding}
    \end{adjustbox}
\end{table}
\begin{table}[t]
    \centering
    \footnotesize
    \caption{Performance of expected value approach and stochastic program \textbf{including} bidding: Objective value (Obj.), net electricity sold on day-ahead market (sales-purchase) (El.sales), Income from electricity market (Income), Heat production (Heat) and cost per MWh. All values are expected values. The VSS is the difference between objective values and VSS[\%] gives the improvement when using the stochastic program.  }
    \begin{adjustbox}{width=\textwidth}
    \begin{tabular}{lrrrrr|rrrrr|rr}
     \toprule
        ~ & \multicolumn{5}{c}{Expected value approach}& \multicolumn{5}{c}{Stochastic program} & ~ & ~ \\ \cmidrule(rl){2-6} \cmidrule(rl){7-11}
        Case & Obj. & El.sales & Income & Heat & Cost/ & Obj. & El.sales & Income & Heat & Cost/ & VSS & VSS\\
          & [EUR]  & [MWh$_\text{e}$] & [EUR] & [MWh$_\text{h}$]& MWh$_\text{h}$ & [EUR]  & [MWh$_\text{e}$] & [EUR] & [MWh$_\text{h}$]& MWh$_\text{h}$ & [EUR]& [\%] \\\midrule
        %\rot{Case} & \rot{Objective value [EUR]} & \rot{El. sales[MWh]} &\rot{Market income [EUR]} & \rot{Heat production [MWh]} & \rot{Cost per MWh-heat} & \rot{Objective value [EUR]} & \rot{El. sales[MWh]} & \rot{Market income [EUR]} & \rot{Heat production [MWh]} & \rot{Cost per MWh-heat} & \rot{VSS [EUR]} & \rot{VSS [\%]} \\ 
        B-01-168 & 144497.6 & 824.7 & 25553.4 & 4956.9 & 29.2 & 106370.7 & 834.4 & 65405.3 & 4957.2 & 21.5 & 38126.9 & 26.4\% \\ 
        B-04-168 & 66262.4 & -33.7 & -28814.4 & 3543.9 & 18.7 & 32010.1 & -3.8 & 11645.2 & 3544.0 & 9.0 & 34252.3 & 51.7\% \\ 
        B-06-168 & 33424.1 & -67.9 & -22584.0 & 2197.0 & 15.2 & 14327.4 & -71.6 & -3494.7 & 2197.8 & 6.5 & 19096.7 & 57.1\% \\ 
        H-04-168 & 288430.3 & 2133.6 & 182284.9 & 8490.5 & 34.0 & 267060.0 & 3768.1 & 307183.5 & 8490.9 & 31.5 & 21370.3 & 7.4\% \\ 
        H-07-168 & 29094.9 & 1258.2 & 149114.0 & 2665.8 & 10.9 & -561.4 & 2347.5 & 267112.2 & 2910.3 & -0.2 & 29656.4 & 101.9\% \\ 
        H-10-168 & -339867.2 & 5037.7 & 914679.2 & 6245.8 & -54.4 & -520133.3 & 9111.8 & 1523616.5 & 10243.8 & -50.8 & 180266.0 & 53.0\% \\ 
        M-05-168 & 16270.5 & 0.0 & 0.0 & 668.5 & 24.3 & 16270.6 & 0.0 & 0.0 & 668.5 & 24.3 & -0.1 & 0.0\% \\ 
        M-08-168 & 6899.4 & 18.2 & 1142.4 & 291.1 & 23.7 & 6862.0 & 28.0 & 1722.4 & 291.1 & 23.6 & 37.4 & 0.5\% \\ 
        M-12-168 & 33603.0 & 85.1 & 4593.1 & 1247.7 & 26.9 & 32799.8 & 180.2 & 9008.5 & 1247.6 & 26.3 & 803.2 & 2.4\% \\  \bottomrule
    \end{tabular}
    \end{adjustbox}
     \label{tab:vss-bidding}
\end{table}

 In the absence of power market bidding (Table \ref{tab:vss-nobidding}), using the stochastic program is particularly beneficial in  the Brønderslev and Hillerød DH systems. In Middelfart, the performance of the stochastic and expected value models are similar. The stochastic program achieves lower total costs by utilizing cheaper units as can be seen by the lower cost per MWh produced heat.  Although the expected value approach achieves higher income from the electricity market in some cases, the cost per MWh heat are always higher or equal. Including scenarios in the stochastic program enables the model to schedule the first-stage units such that total costs across all scenarios are decreased in comparison to the expected value program. The expected value approach optimizes only for the expected scenario and therefore, the commitment of units can be disadvantageous for other scenarios. In Middelfart, the interaction with the electricity market is low, since the cheapest units (wood chip boiler and wood pellet boiler) do not depend on the market. The uncertainty can be handled successfully using the flexibility of the system itself (heat storage and flexible, market-independent boilers), so the costs are the same for both approaches.
 
 If we drop the assumption that the units can obtain the day-ahead market prices without bidding, the benefits of the stochastic program increase further. Table \ref{tab:vss-bidding} shows the results when the bidding behaviour is modelled using bidding curves. Then, even for Middelfart, some benefits can be achieved. The modelling of bidding curves relies on the scenarios for electricity prices, i.e., for each  price scenario a single price and quantity are determined. In the expected value setting, only one scenario is present, hence only one bid is given. The distinction of several production levels  enables us to achieve more won bids and,  consequently, higher profits from the electricity market which reduces the overall costs. In case H-10-168, the electricity prices are so favorable that the stochastic program produces more heat than needed to create additional income from the market.
 
 When comparing Table \ref{tab:vss-nobidding} and Table \ref{tab:vss-bidding}, we see that constructing bidding curves leads to the stochastic program achieving lower in-sample costs in all cases. The bidding curves allow us to distinguish different production levels, while the operational model without bidding curves has to determine one level of operation for all units with commitment decisions and periods.
 
 In both tables, seasonal variations can be observed. The cost per MWh are lower in summer where the heat demand is lower and thus cheaper units are enough to cover the demand. Depending on the cost structure in the respective DH system, this leads to a different model behaviour across seasons: In Middelfart and Br{\o}nderslev, market-independent units tend to yield the lowest cost and thus, trading volumes and per-MWh cost of heat are higher during winter, when market-dependent units need to be active to cover the higher heat demand. In Hiller{\o}d, on the other hand, the market-dependent CHP unit tends to be a low-cost heat generator throughout the entire year. In October 2021, electricity prices rise and therefore, the heat cost per MWh can be lower due to extended electricity sales at higher prices.

  \begin{figure}
 	\centering
 	\begin{subfigure}[b]{0.49\columnwidth}
 		\centering
 		\includegraphics[width=0.95\columnwidth]{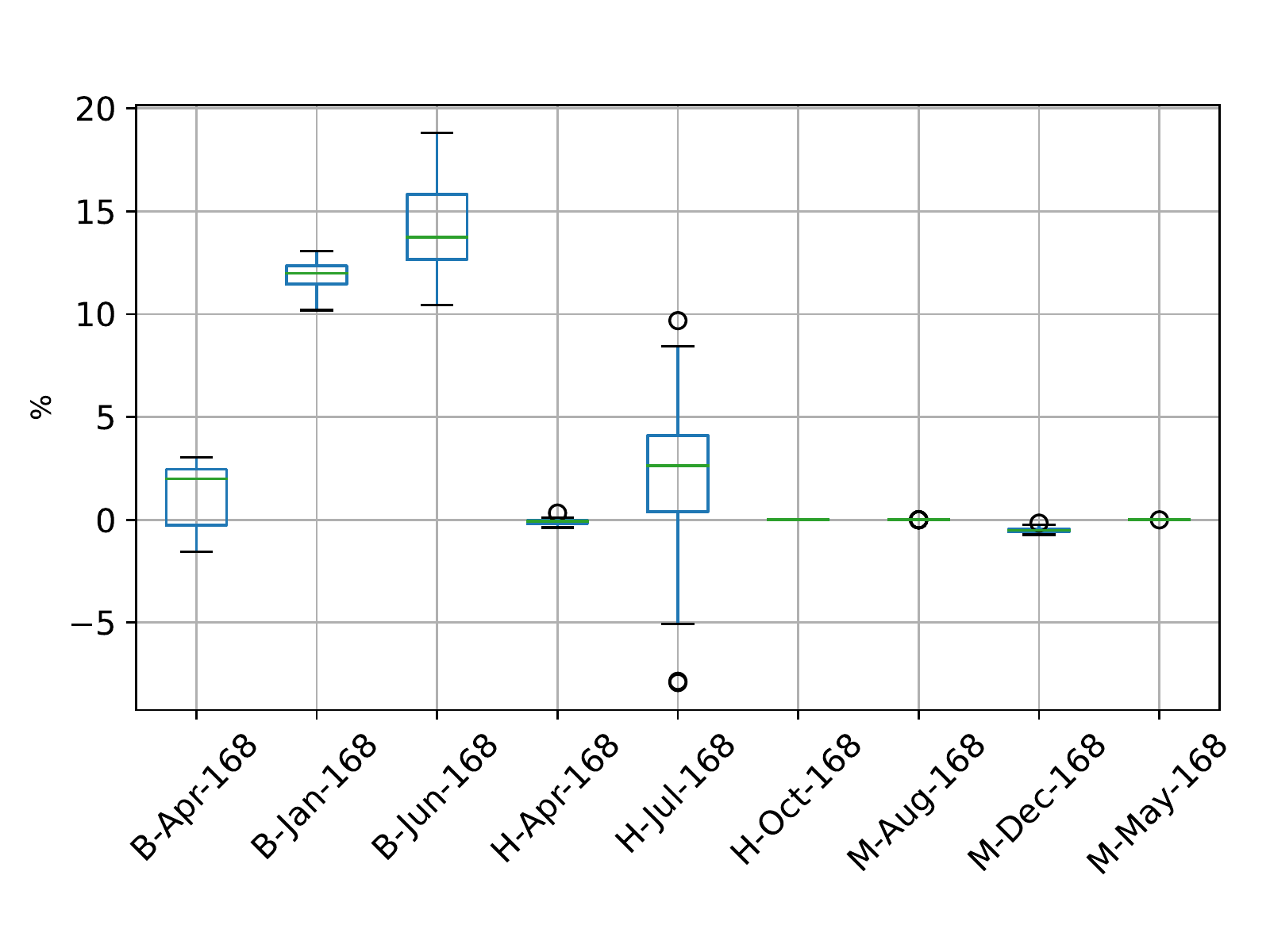}
 		\caption{Without bidding}
 		\label{fig:oos-no-bid}
 	\end{subfigure}
 	\hfill
 	\begin{subfigure}[b]{0.49\columnwidth}
 		\centering
 		\includegraphics[width=0.95\columnwidth]{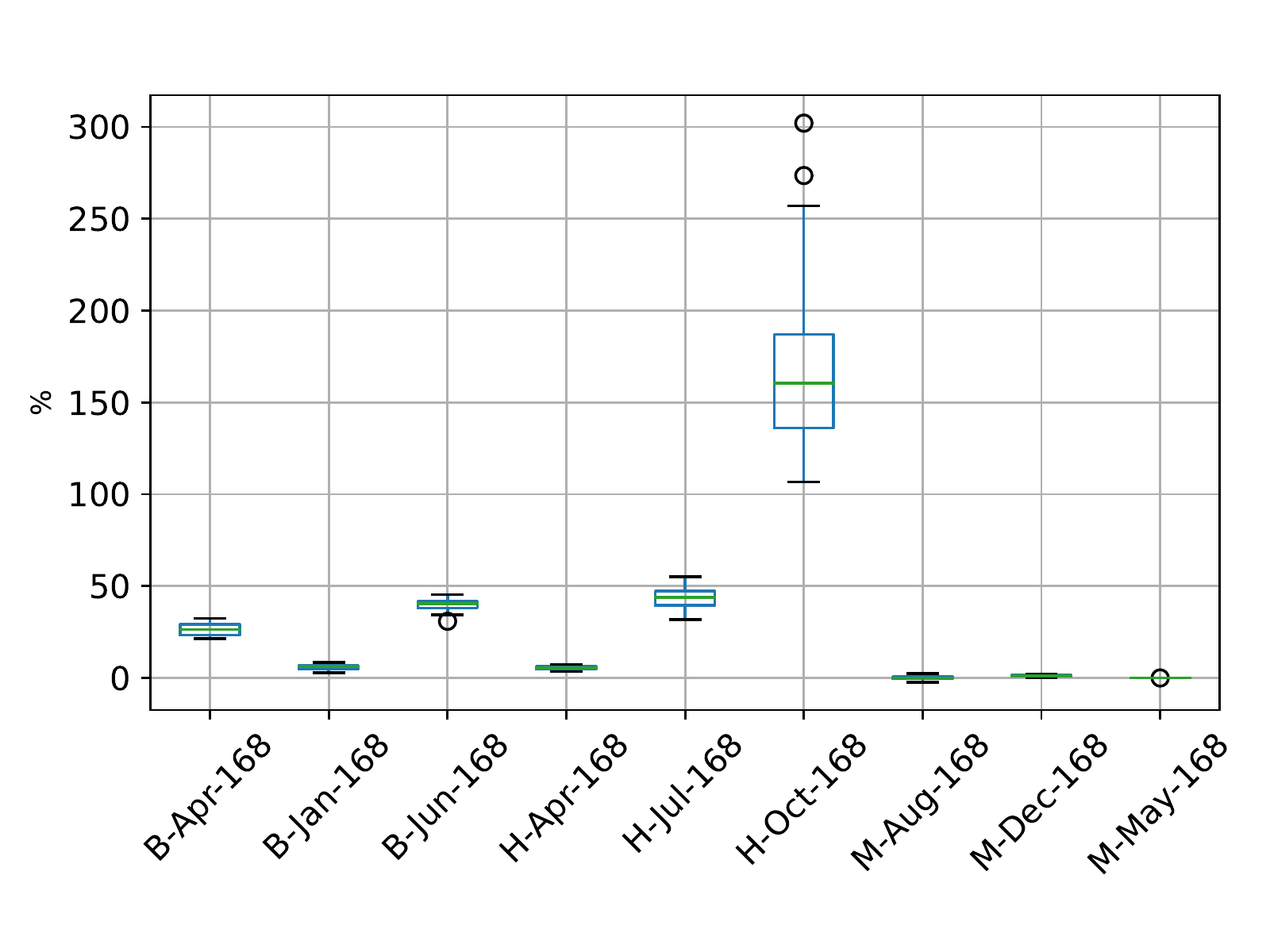}
 		\caption{With bidding}
 		\label{fig:oos-bid}
 	\end{subfigure}
 	\caption{Out-of-sample performance for set 1. Values are the difference between exp. value approach and stochastic programming in \%. Plots contain the values of the 30 sampled realizations of uncertainty.}
 	\label{fig:oos}
 \end{figure}
 
 \begin{figure}
 	\centering
 	\begin{subfigure}[b]{0.49\columnwidth}
 		\centering
 		\includegraphics[width=0.95\columnwidth]{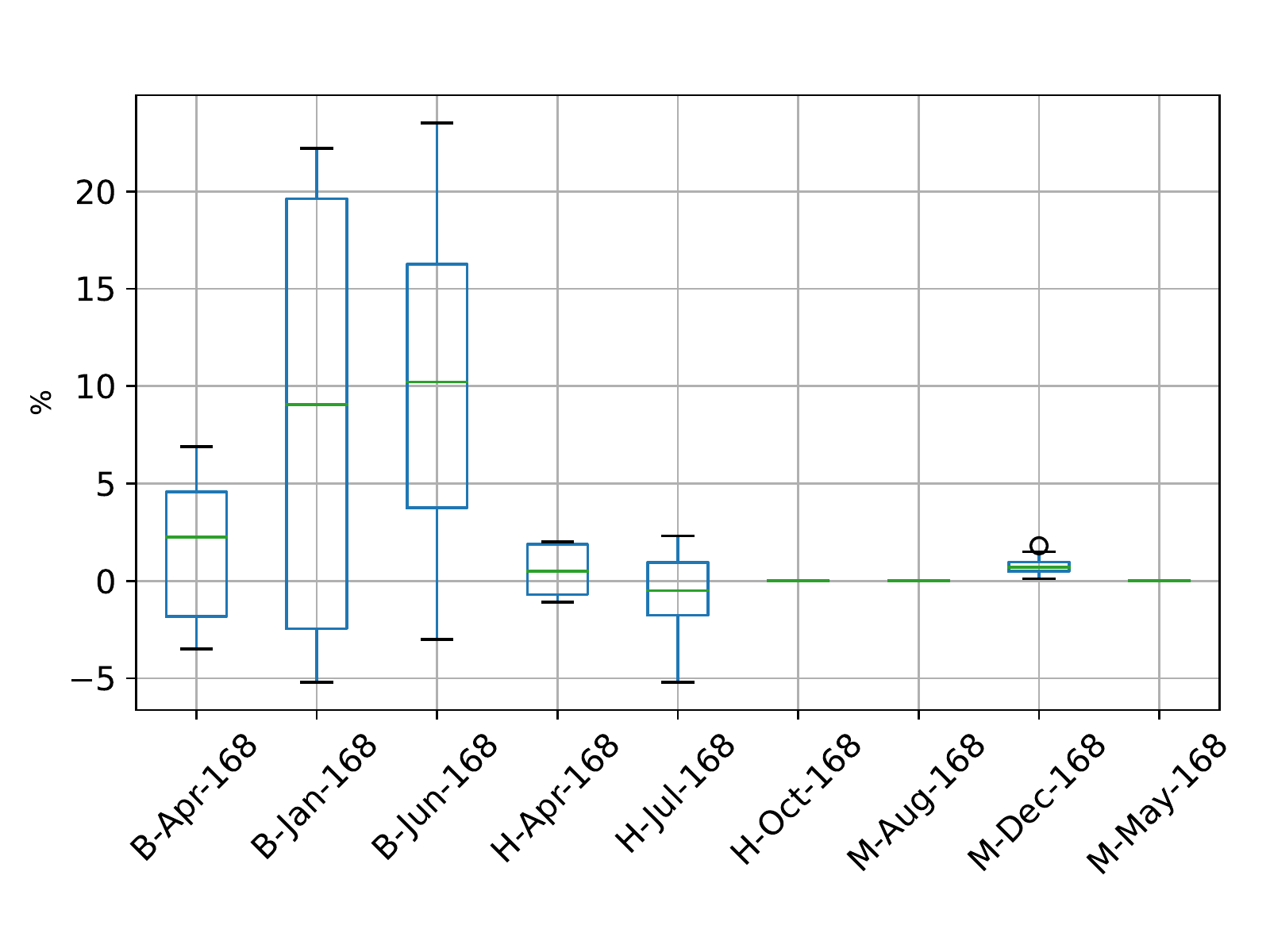}
 		\caption{Without bidding}
 		\label{fig:oos2-no-bid}
 	\end{subfigure}
 	\hfill
 	\begin{subfigure}[b]{0.49\columnwidth}
 		\centering
 		\includegraphics[width=0.95\columnwidth]{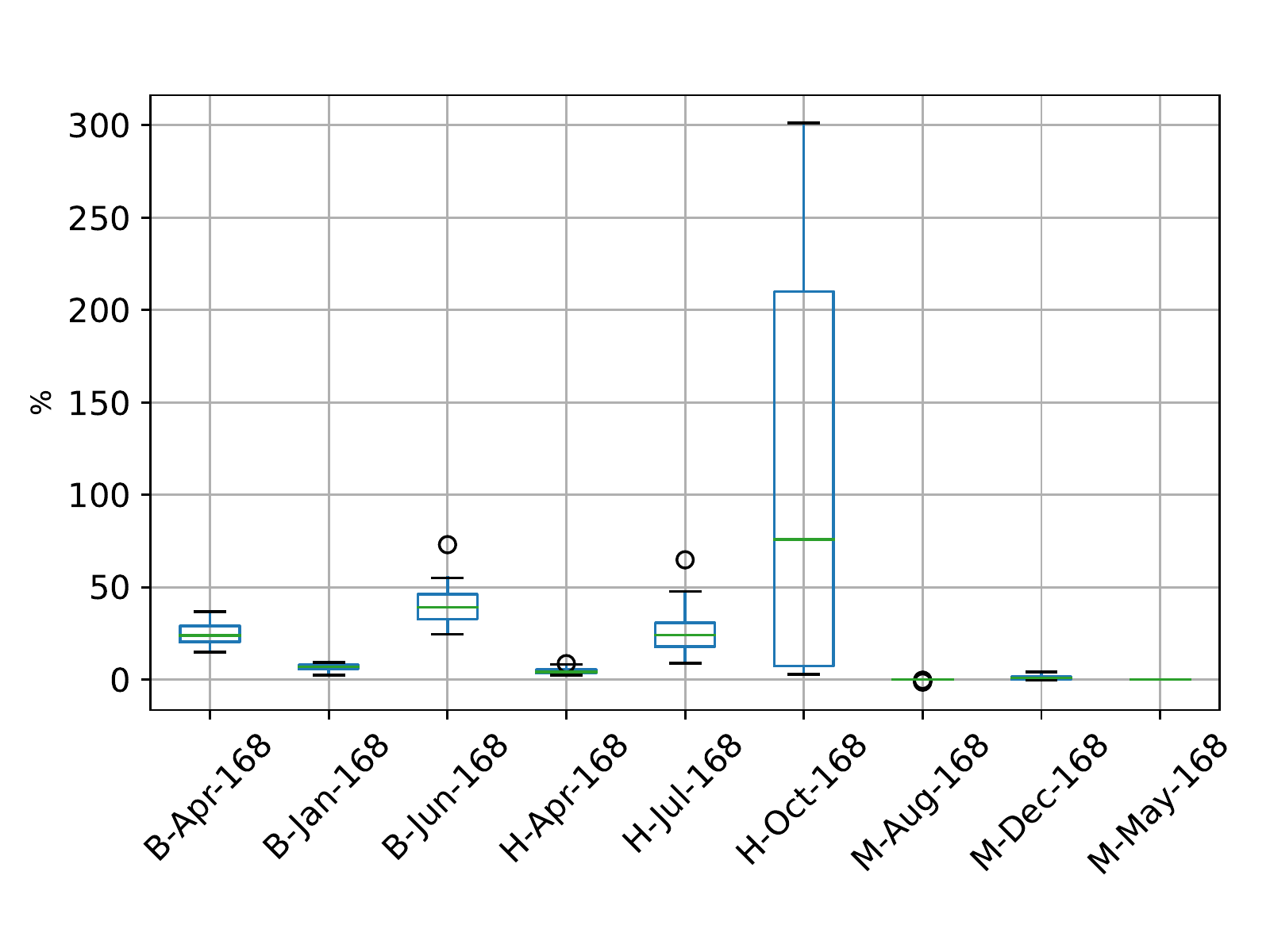}
 		\caption{With bidding}
 		\label{fig:oos2-bid}
 	\end{subfigure}
 	\caption{Out-of-sample performance for set 2. Values are the difference between exp. value approach and stochastic programming in \%. The plots contain the values of the 30 sampled realizations of uncertainty.}
 	\label{fig:oos2}
 \end{figure}
 
%  In colder months, the costs and also the trading volume (El. sales) increase. Due to the higher heat demand, the market-dependent units need to be used to produce heat and the models schedule the production in hours with favorable expected electricity prices.
 
 Since the VSS only uses scenarios that were already considered in the model, we additionally perform an out-of-sample evaluation. Here, the first-stage decisions are evaluated using new samples that were not used in the initial optimization to account for robust performance in unseen cases. For that purpose, we generate two sample sets:
 \begin{enumerate}

     \item In the first set, we sample the uncertain parameters from a triangular distribution for each hour of the week. The parameters of the triangular distribution are specific for each hour to account for different patters during the day and week. We estimate the parameters using the values from the same hour during the week in the three historic weeks used in the scenarios. We deduct/add additional 5\% from the lowest/to the highest value to allow also for slightly lower/higher values in the sampling.
     The set contains 30 weekly samples per case.
     \item In the second set, we use block bootstrapping with historic data from the two weeks prior and the two weeks after the three weeks used in the scenarios, i.e., four weeks in total. Afterwards, we split the data into blocks of four hours to capture some temporal dependency. To create new cases, we sample six blocks of four hours per day. Each block is chosen randomly from the available blocks, but the selection is limited to the same time of the day and distinguishes between weekday and weekend. The set contains 30 weekly samples per case, 15 cases sampled from the first two weeks and 15 cases from the latter two weeks.
 \end{enumerate}
 
 The comparison of the expected value model and stochastic model for sets 1 and 2 is presented in Figures \ref{fig:oos} and \ref{fig:oos2}, respectively. The shown data points are the savings of the stochastic program compared to the expected value model (in \%), i.e., positive values mean that the stochastic programming solution outperformed the expected value solution.

 The results confirm the observations from the VSS results, i.e, in particular Brønderslev benefits from using a stochastic program and the benefits increase when considering bidding curves. When no bidding is considered, the expected value approach outperforms the stochastic program in some cases, but on average the stochastic program leads to better results. When bidding is considered, the stochastic program clearly outperforms the expected value approach. Based on the out-of-sample evaluation, we can conclude that the stochastic program outperforms an expected value approach. In DH systems with less interaction on the market and no uncertain production (such as Middelfart), a deterministic model and a simple point forecast might be sufficient, while in other cases the stochastic program should be utilized. Note that the scenario generation used in this paper is very simple. Using scenarios created from proper probabilistic forecasting techniques will increase the benefit of using the stochastic program in many cases.

\subsection{Evaluation on real data with rolling horizon}\label{sec:rc}

Next, we evaluate the operational planning in a rolling horizon setting for the same 9 scenarios as previously. Afterwards, the first-stage solutions are evaluated on the realization of the uncertain data. The rolling horizon approach applies a sliding window with length of 168h and we move the window by 24h after each iteration, i.e., the non-anticipativity period for the first-stage decisions is 24h. We consider 168h in the model to account for storage behaviour (as shown in \cite{hurb}). In the evaluation, the storage levels and unit statuses at the end of each day are determined using the observed data and applied as input for the next optimization. Thus, we mimic daily planning in practice. In total, we optimize for 2 weeks, i.e., 14 iterations, and the planning horizon reduces from 168h in a receding manner when approaching the end of the two weeks.

Table \ref{tab:rh} shows the total cost when applying the optimization results to the real observations over the 14 days. If we do not consider bidding curves, the stochastic programming solution can reduce the cost slightly in most cases. Only in one case (B-04-168), the expected value approach performs slightly better than the stochastic program for this specific realization of uncertain parameters. This can be explained by the runtime limitations, as shown in the next section. When considering bidding, the stochastic program reduces cost in all cases and can save up to 42.1\% of the cost in 14 days. The stochastic model performs particularly well for Brønderslev and Hillerød as already concluded in the previous section.

\begin{table}[t]
	\footnotesize
	\caption{Rolling horizon performance on real data. Difference between exp. value approach and stochastic programming evaluated over the a period of 14 days with a planning horizon 168h and rolling window of 24h. }
	\centering
	\begin{adjustbox}{width=0.9\columnwidth}
		\begin{tabular}{lrrrrrrrrr}\toprule
			Case & \multicolumn{4}{c}{No bidding} & \multicolumn{4}{c}{Bidding}\\\cmidrule(r){2-5}\cmidrule(r){6-9}
			& Exp. & Sto. & $\Delta$ & [\%] & Exp. & Sto. & $\Delta$ & [\%] \\ \midrule
			B-01-168 & 257749.88 & 256180.52 & 1569.36 & 0.6\% & 321250.27 & 288560.77 & 32689.50 & 10.2\% \\ 
			B-04-168 & 41257.64 & 41446.22 & -188.58 & -0.5\% & 80486.26 & 55179.63 & 25306.63 & 31.4\% \\ 
			B-06-168 & 21740.36 & 15866.63 & 5873.72 & 27.0\% & 35248.22 & 20422.54 & 14825.68 & 42.1\% \\ 
			H-04-168 & 605527.16 & 603757.54 & 1769.63 & 0.3\% & 531522.78 & 512361.63 & 19161.15 & 3.6\% \\ 
			H-07-168 & -103030.15 & -107056.46 & 4026.31 & 3.9\% & -99279.96 & -127250.59 & 27970.63 & 28.2\% \\ 
			H-10-168 & -1187638.56 & -1187638.56 & 0.00 & 0.0\% & -1199759.10 & -1379174.08 & 179414.98 & 15.0\% \\ 
			M-05-168 & 26082.61 & 26069.20 & 13.41 & 0.1\% & 26082.61 & 26069.20 & 13.41 & 0.1\% \\ 
			M-08-168 & 15532.20 & 15512.30 & 19.90 & 0.1\% & 14687.92 & 14688.02 & -0.10 & 0.0\% \\ 
			M-12-168 & 73741.12 & 72524.16 & 1216.96 & 1.7\% & 72837.97 & 72581.06 & 256.91 & 0.4\% \\    \bottomrule
		\end{tabular}
	\end{adjustbox}
	\label{tab:rh}
\end{table}

\begin{figure}[t]
	\centering
	\begin{subfigure}[b]{0.49\columnwidth}
		\centering
		\includegraphics[width=1.1\columnwidth]{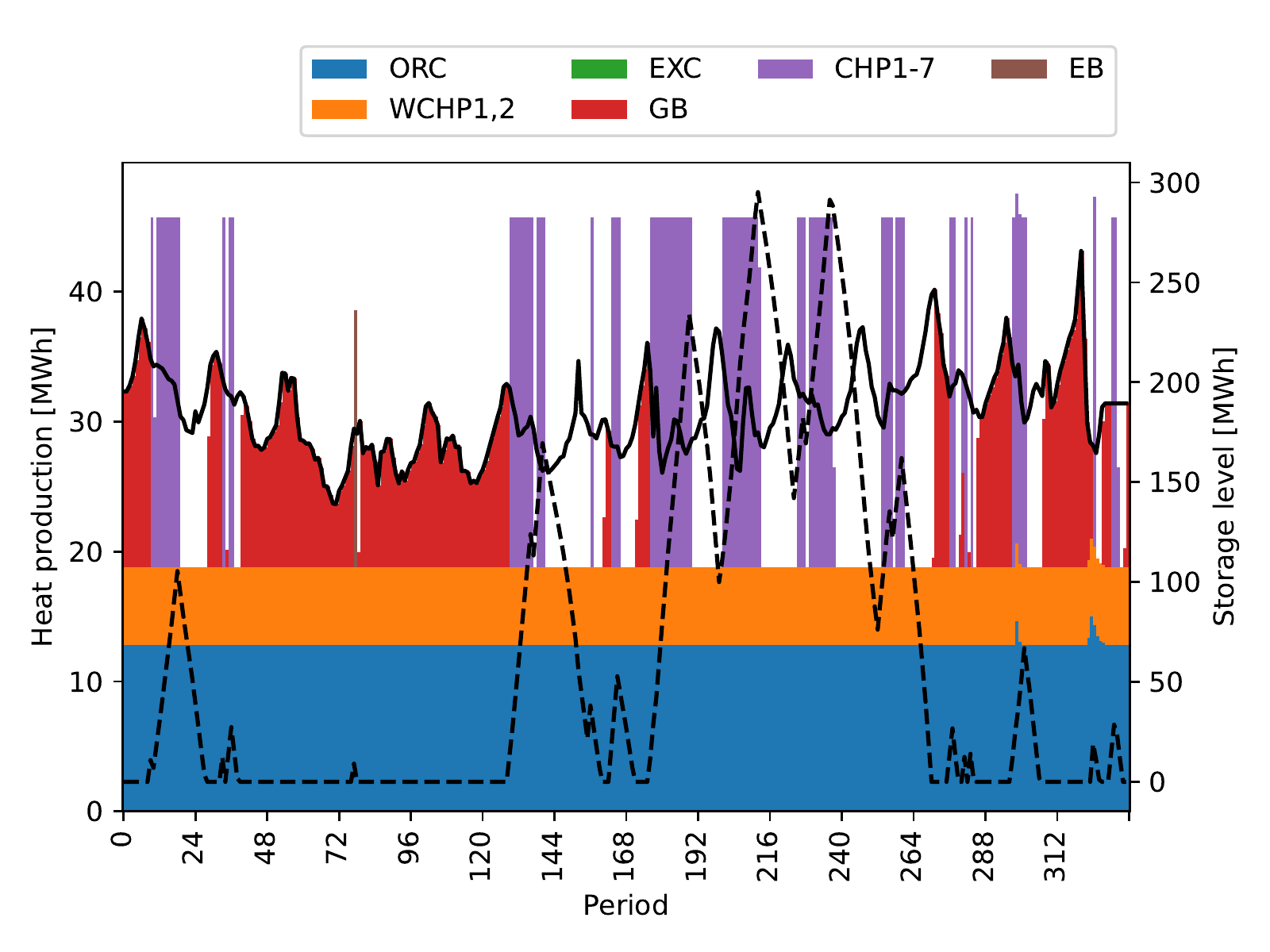}
		\caption{Heat production}
		\label{fig:b-heat}
	\end{subfigure}
	\hfill
	\begin{subfigure}[b]{0.49\columnwidth}
		\centering
		\includegraphics[width=1.1\columnwidth]{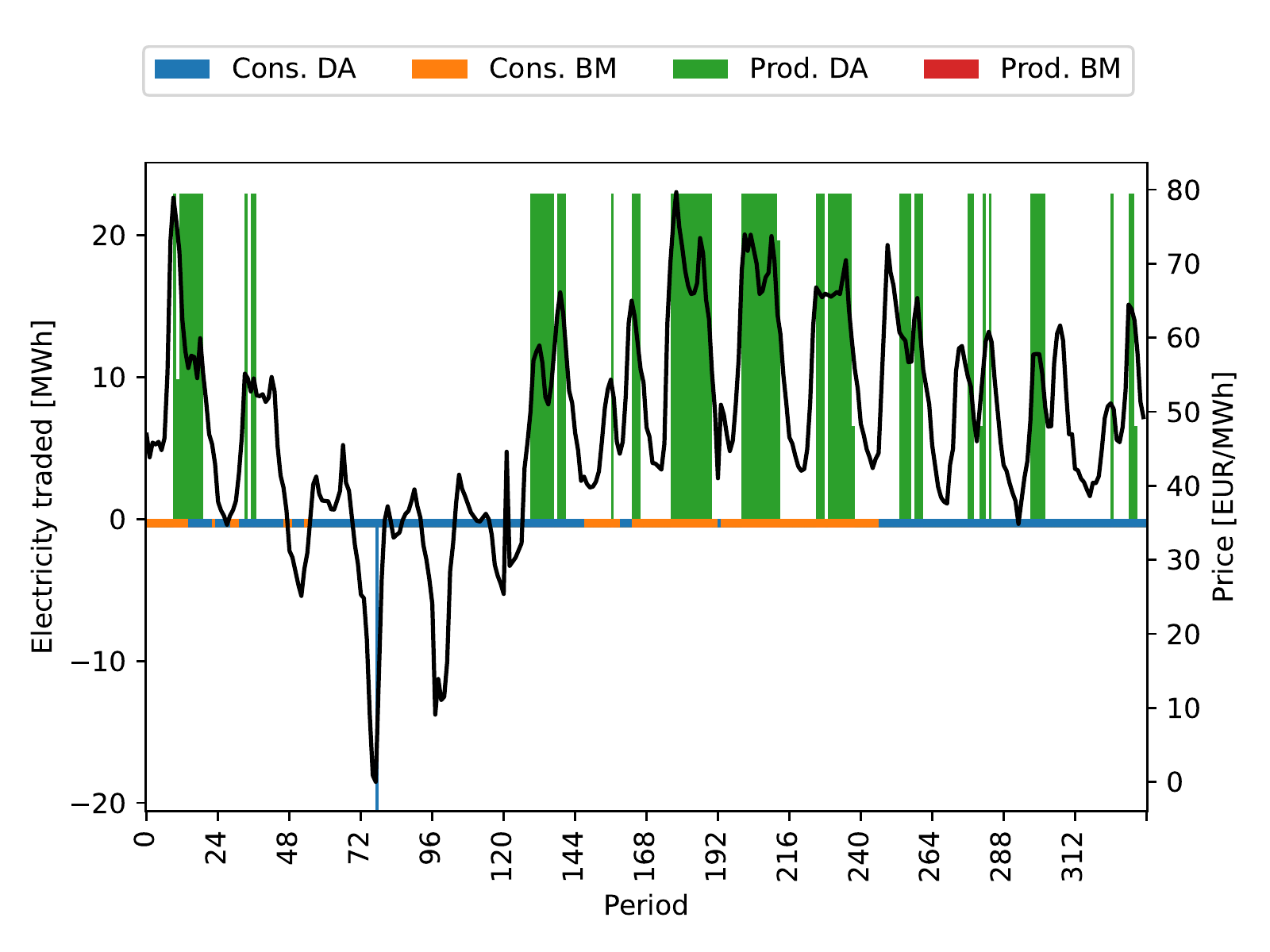}
		\caption{Electricity trading}
		\label{fig:b-el}
	\end{subfigure}
	\caption{Optimized production for case B-01-168 with daily rolling horizon and 168h planning horizon.}
	\label{fig:b-01-168}
\end{figure}
\begin{figure}[t]
	\centering
	\begin{subfigure}[b]{0.49\columnwidth}
		\centering
		\includegraphics[width=1.1\columnwidth]{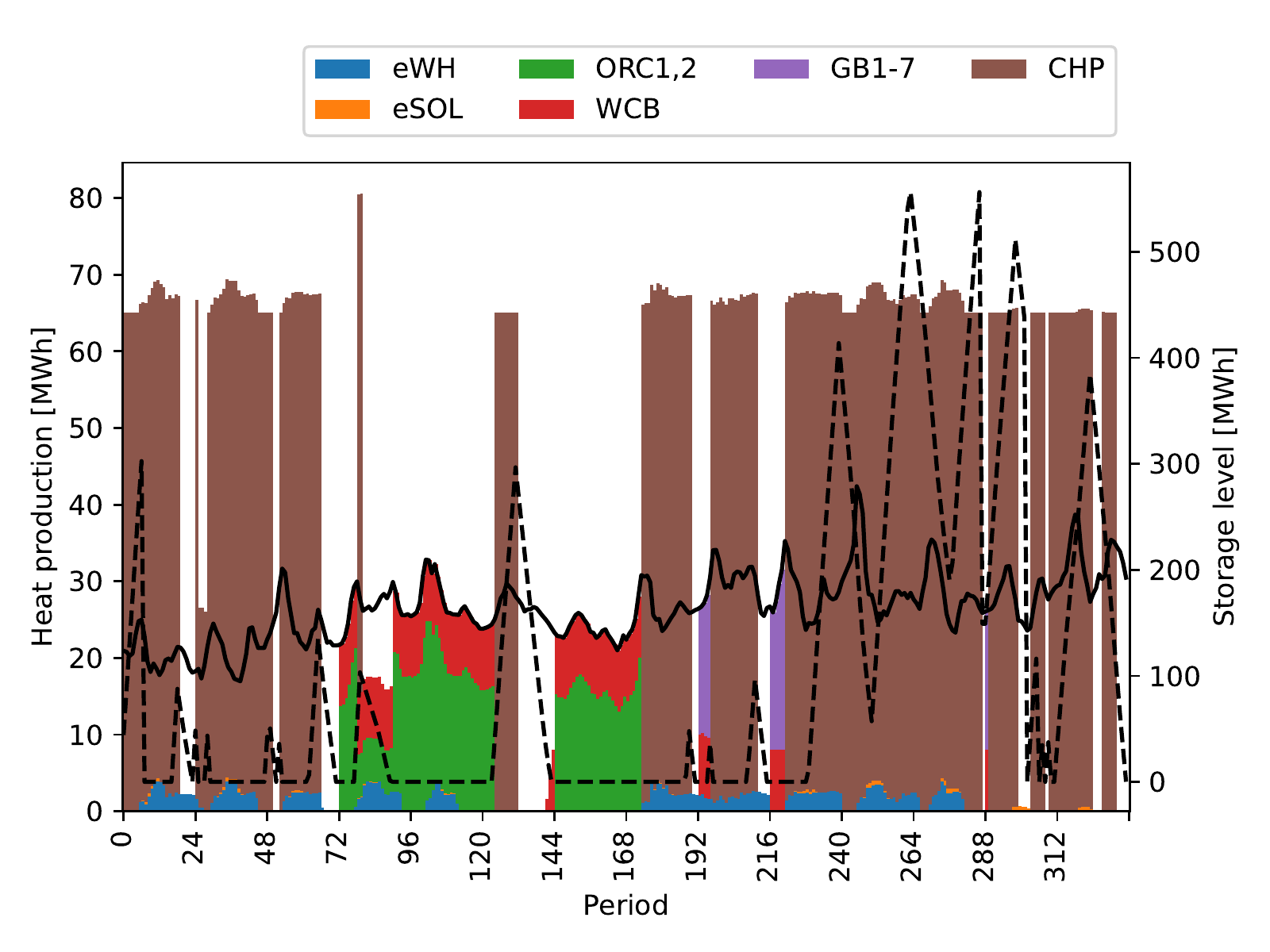}
		\caption{Heat production}
		\label{fig:h-heat}
	\end{subfigure}
	\hfill
	\begin{subfigure}[b]{0.49\columnwidth}
		\centering
		\includegraphics[width=1.1\columnwidth]{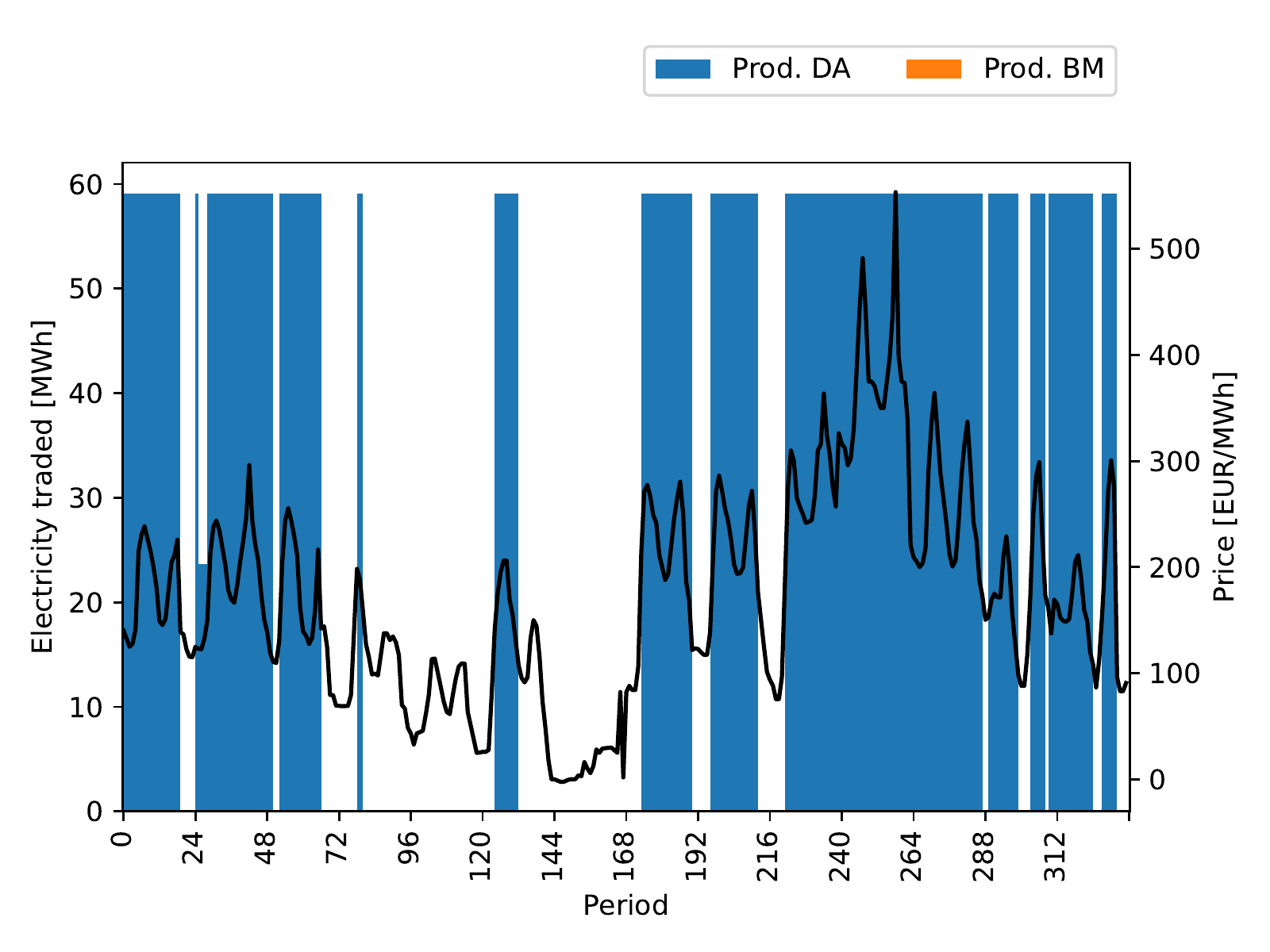}
		\caption{Electricity trading}
		\label{fig:h-el}
	\end{subfigure}
	\caption{Optimized production for case H-10-168 with daily rolling horizon and 168h planning horizon.}
	\label{fig:h-04-168}
\end{figure}
\begin{figure}[ht]
	\centering
	\begin{subfigure}[b]{0.49\columnwidth}
		\centering
		\includegraphics[width=1.1\columnwidth]{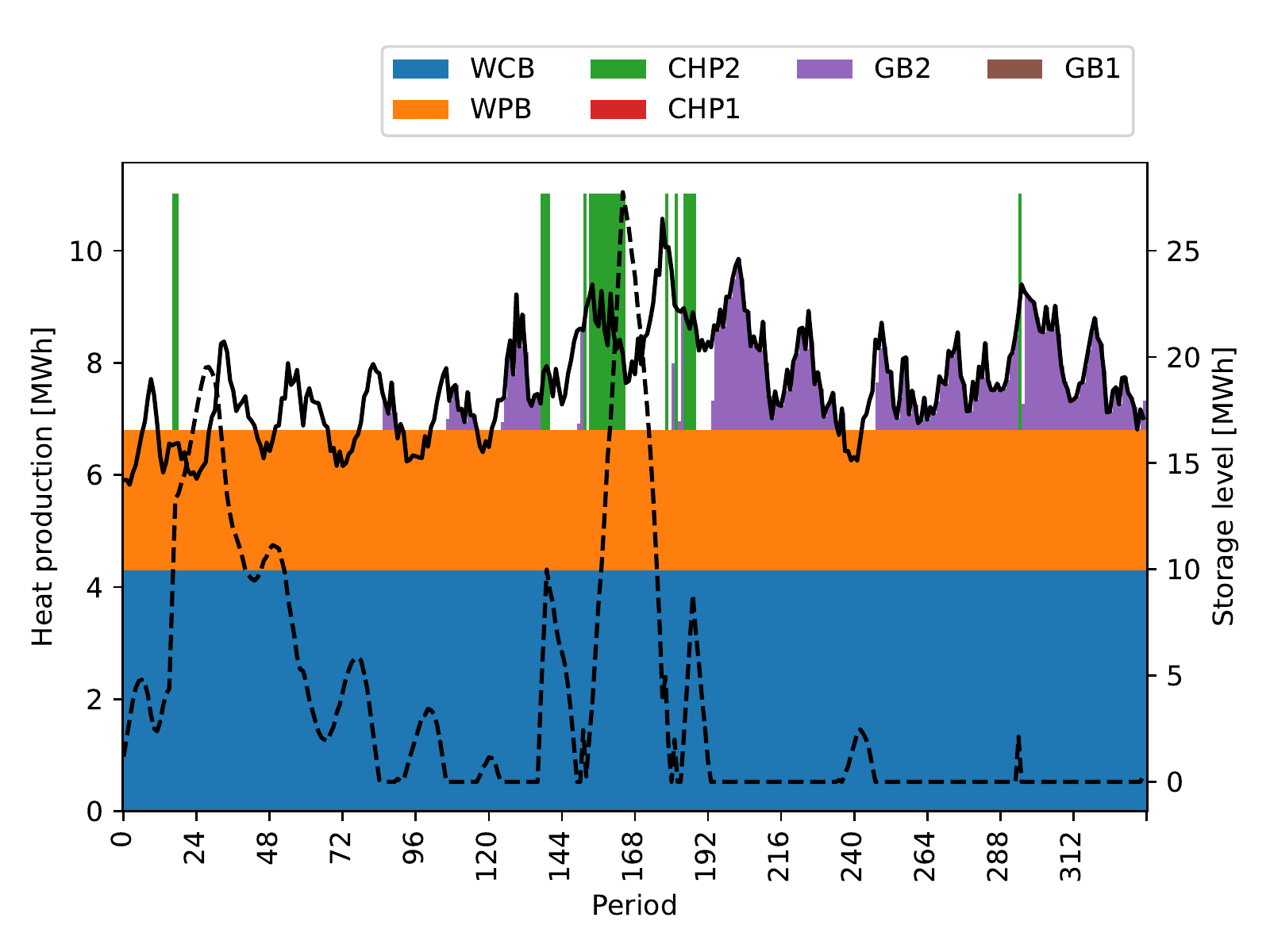}
		\caption{Heat production}
		\label{fig:m-heat}
	\end{subfigure}
	\hfill
	\begin{subfigure}[b]{0.49\columnwidth}
		\centering
		\includegraphics[width=1.1\columnwidth]{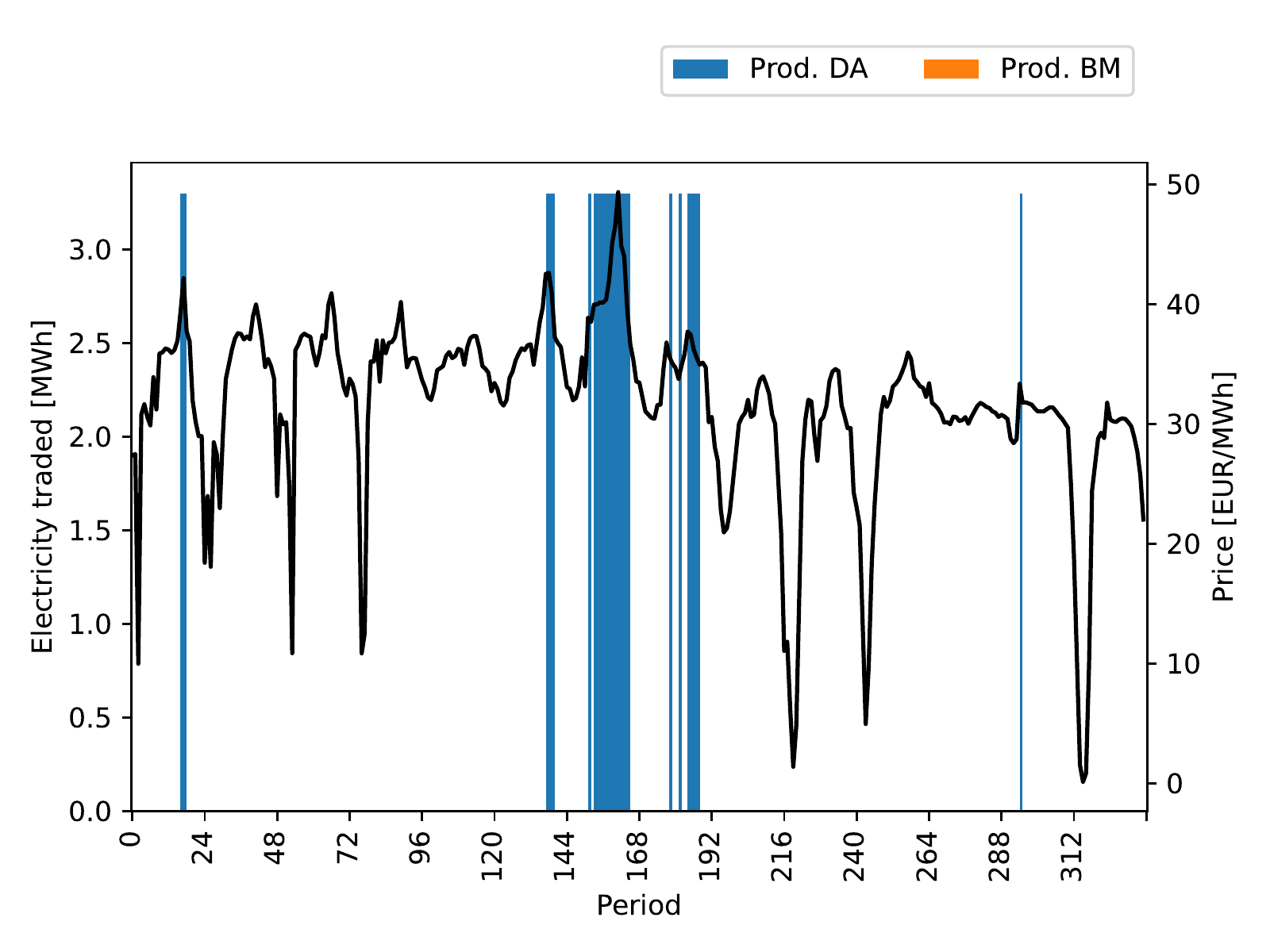}
		\caption{Electricity trading}
		\label{fig:m-el}
	\end{subfigure}
	\caption{Optimized production for case M-12-168 with daily rolling horizon and 168h planning horizon.}
	\label{fig:m-12-168}
\end{figure}

Figures \ref{fig:b-01-168}, \ref{fig:h-04-168} and \ref{fig:m-12-168} show the heat production and electricity trading for 14 days for one case in Brønderslev, Hillerød and Middelfart, respectively. The results are based on the model using bidding curves.
In Brønderslev (Fig. \ref{fig:b-01-168}), the heat pumps (WCHP) with the ORC unit provide the base load for heat production. The remaining heat demand is covered either by the gas boiler (GB) or the CHP (CHP1-7) units. The CHP units fill the storage unit in hours where electricity is sold and storage outflow covers the heat demand in the succeeding hours. On the third day, the electric boiler (EB) wins a bid, since the electricity price is very low (see Fig \ref{fig:b-el}). The CHP units win bids and produce in hours with high electricity prices. The electricity needed for the heat pumps is either won as bids on the day-ahead market (Cons. DA) or imbalances are created (Cons. BM) (see Fig \ref{fig:b-el}) since the combination of heat pumps with the ORC is very cost-efficient.

In Hillerød (Fig. \ref{fig:h-04-168}), the CHP unit also exploits high-price hours to produce heat and fill the storage. The remaining heat demand is either covered by the wood chip boiler (WCB), waste heat (eWH), solar heat (eSOL) or the ORC unit (ORC1,2). The gas boilers are used only in exceptional cases due to their high operational cost.

In Middelfart (Fig. \ref{fig:m-12-168}), the CHP1 and GB1 units are never used, as they are too expensive. Only CHP2 wins production bids in some hours with high prices. The main heat production in Middelfart is achieved using the wood chip boiler (WCB) and wood pellet boiler (WPB). When no CHP bids are successful, peak demand is covered by the gas boiler (GB).

\subsection{Runtimes}

\begin{figure}[ht!]
    \centering
      \begin{subfigure}[b]{0.49\columnwidth}
         \centering
         \includegraphics[width=0.95\columnwidth]{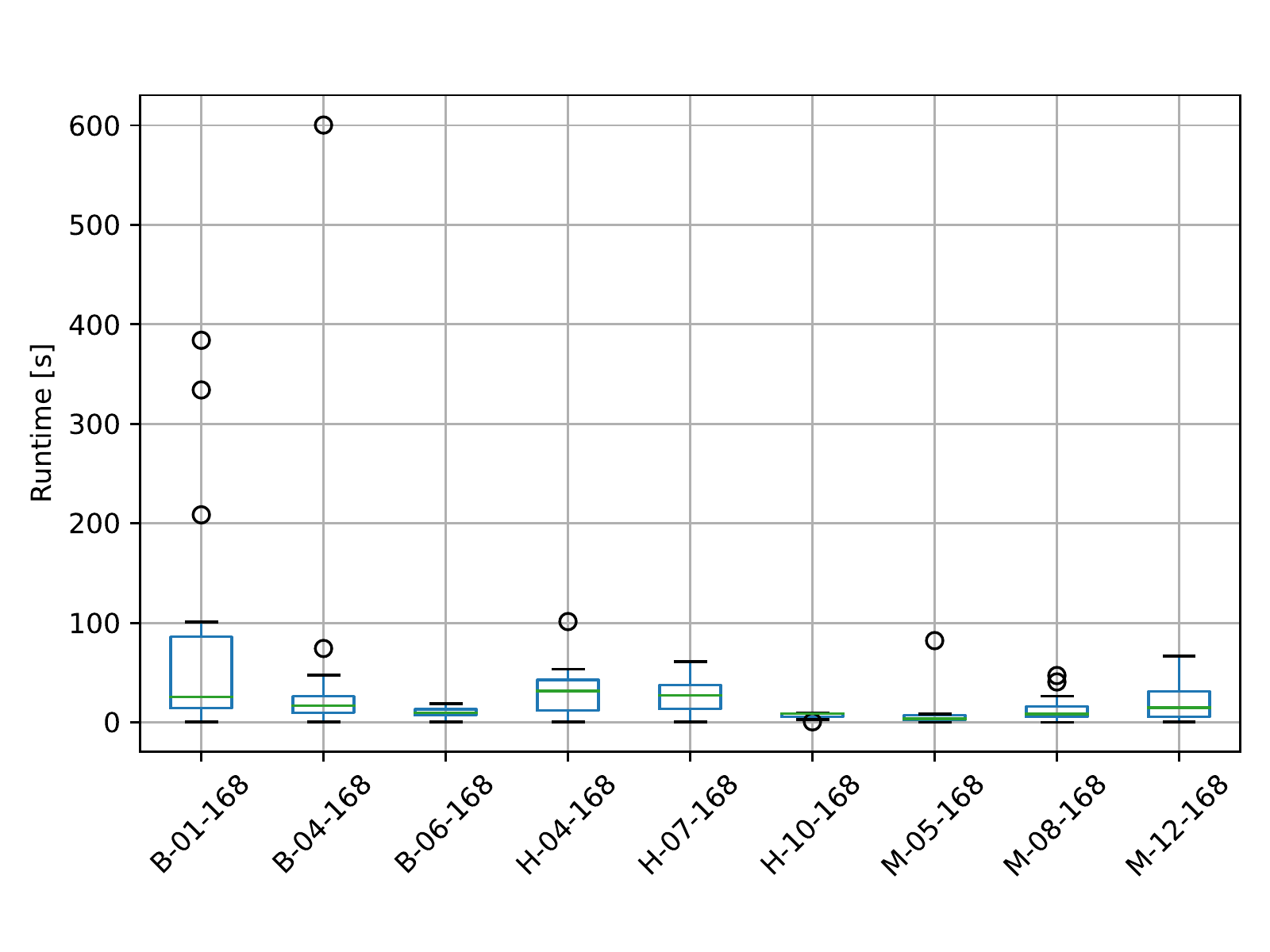}
         \caption{No bidding}
         \label{fig:runtime-nobidding}
     \end{subfigure}
     \hfill
     \begin{subfigure}[b]{0.49\columnwidth}
         \centering
         \includegraphics[width=0.95\columnwidth]{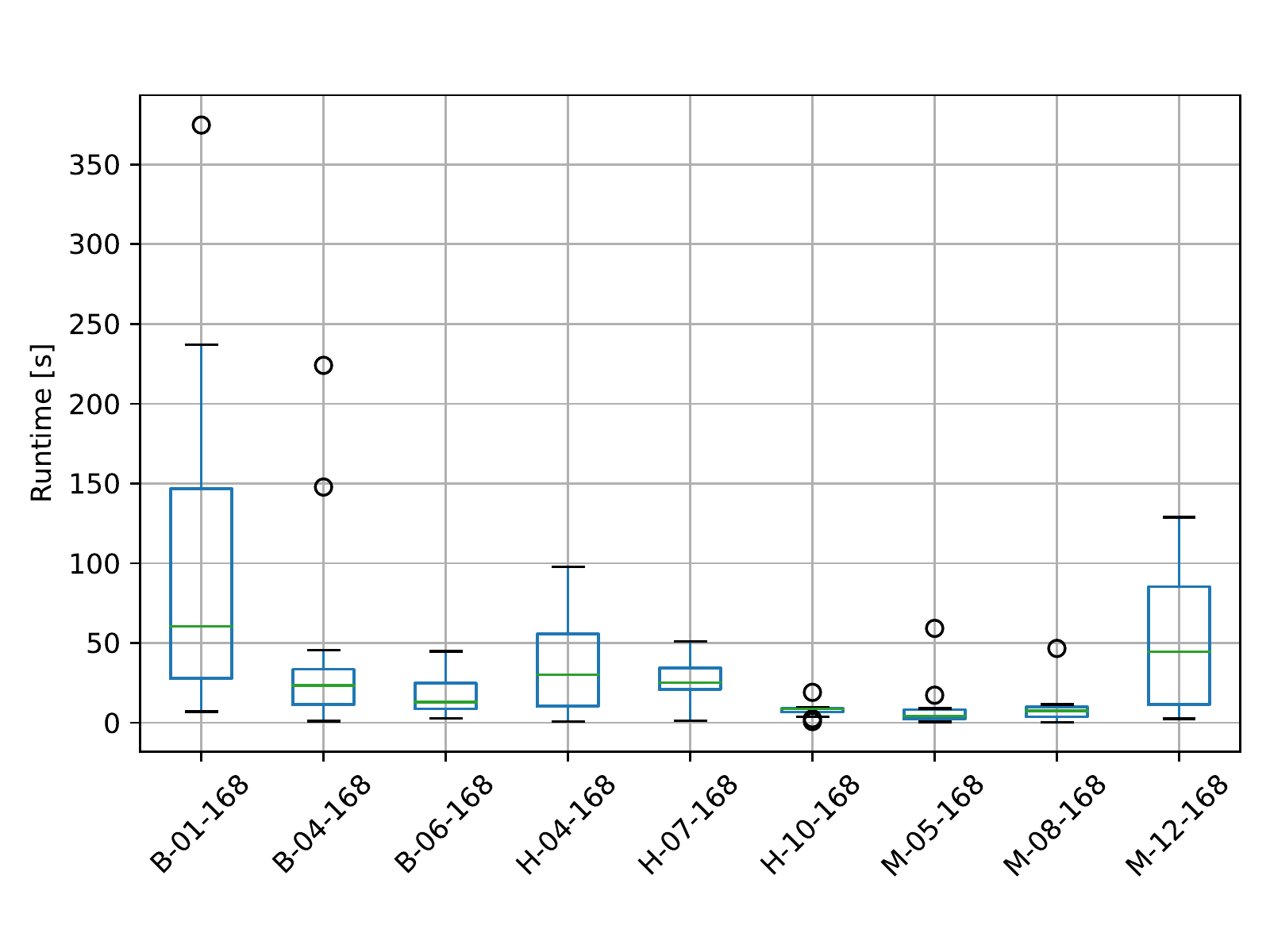}
         \caption{Bidding}
         \label{fig:runtime-bidding}
     \end{subfigure}
    \caption{Computational runtime for each case and all 14 iterations of the rolling horizon. Timeout is 600 sec.}
    \label{fig:runtime}
\end{figure}

The runtimes for solving the model in each iteration of the receding horizon algorithm presented in the previous Section \ref{sec:rc} are given Figure \ref{fig:runtime}. This means each box plot is based on 14 values. The time-out per iteration is 600 seconds and is only hit once in case B-04-168. Brønderslev is by far the most complex network of the three.  For the cases in Hillerød and Middelfart, each model is solved in less than 135 seconds and on average in less than 35 seconds. For the iteration where the time-out of 600 seconds is reached, the remaining gap is 0.0003 (0.03\%) close to the default cutoff of Gurobi (0.0001). Note that B-04-168 is also the case in Table \ref{tab:rh} where the stochastic program did not outperform the expected value approach.

Based on these results, we deem the runtime as short enough to be used in practice, since the optimization is carried out only once a day (for day-ahead market optimization). Even if the generic formulation is used for intra-day optimization, runtimes less than 10 minutes are fast enough.

\subsection{Long-term analysis} \label{subsec:long-term-analysis}
In the last analysis, we apply the model formulation to evaluate the performance of the DH systems in the long term. For this, we use the cases M-03-6936, B-10-6936 and H-02-5808 that include historical data for more than seven months. We solve the entire planning horizon as a deterministic case with the assumption that the day-ahead market can be used without bidding.

\begin{figure}[t]
    \centering
      \begin{subfigure}[b]{0.49\columnwidth}
         \centering
         \includegraphics[width=1\columnwidth]{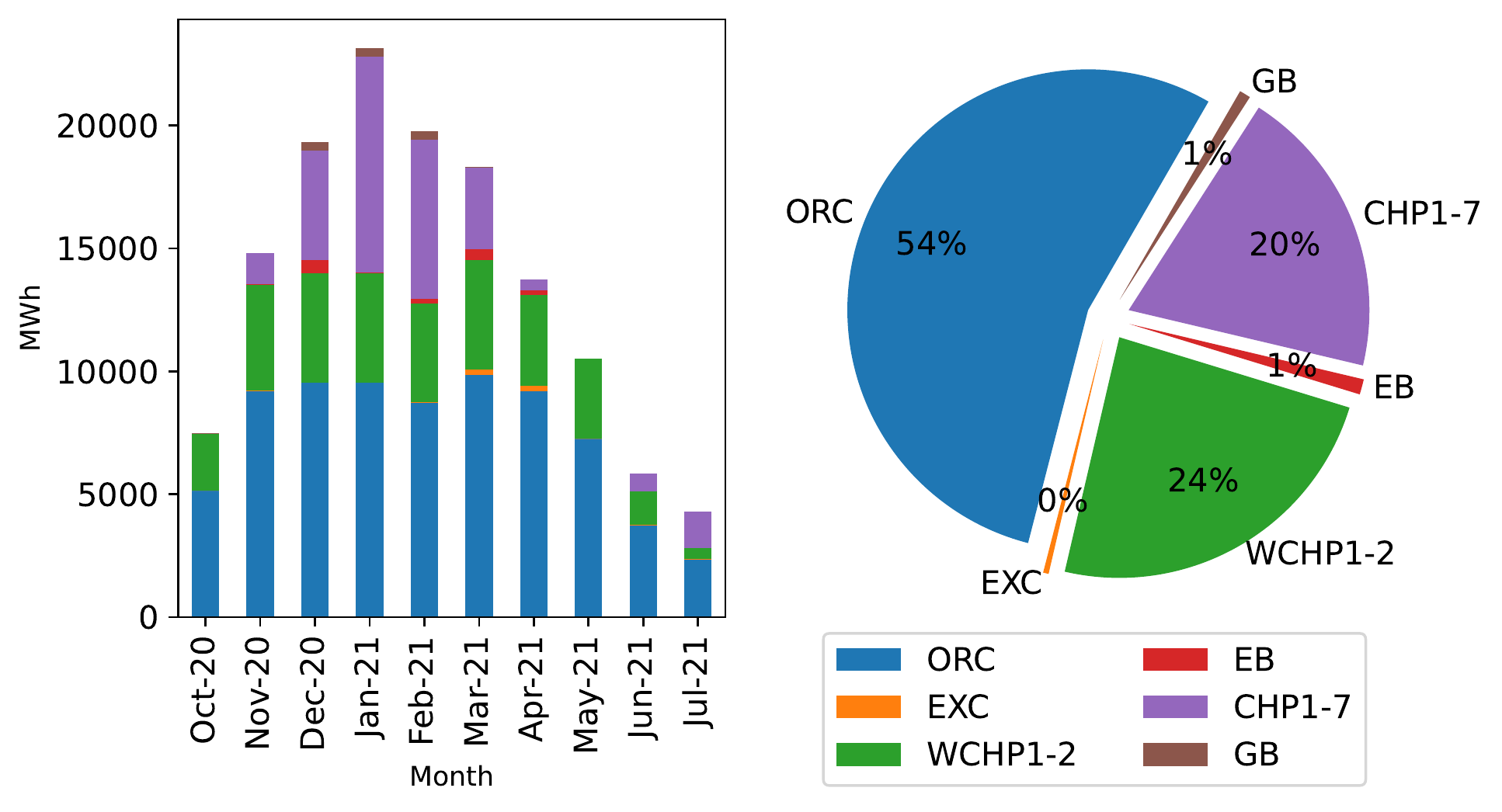}
         \caption{B-10-6936}
         \label{fig:b-total}
     \end{subfigure}
     \hfill
      \begin{subfigure}[b]{0.49\columnwidth}
         \centering
         \includegraphics[width=1\columnwidth]{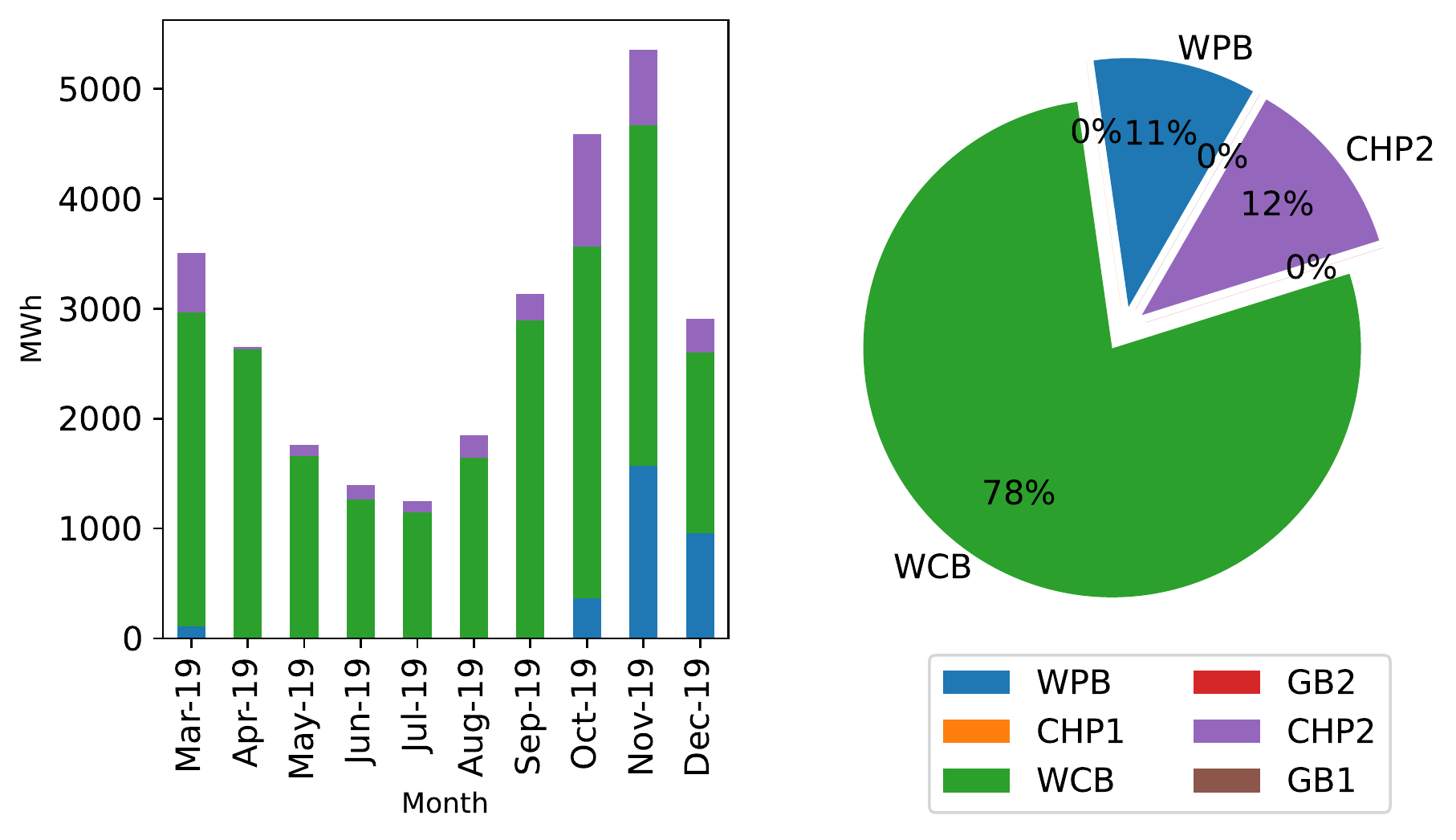}
         \caption{M-03-6936}
         \label{fig:m-total}
     \end{subfigure}
     \begin{subfigure}[b]{0.49\columnwidth}
         \centering
         \includegraphics[width=1\columnwidth]{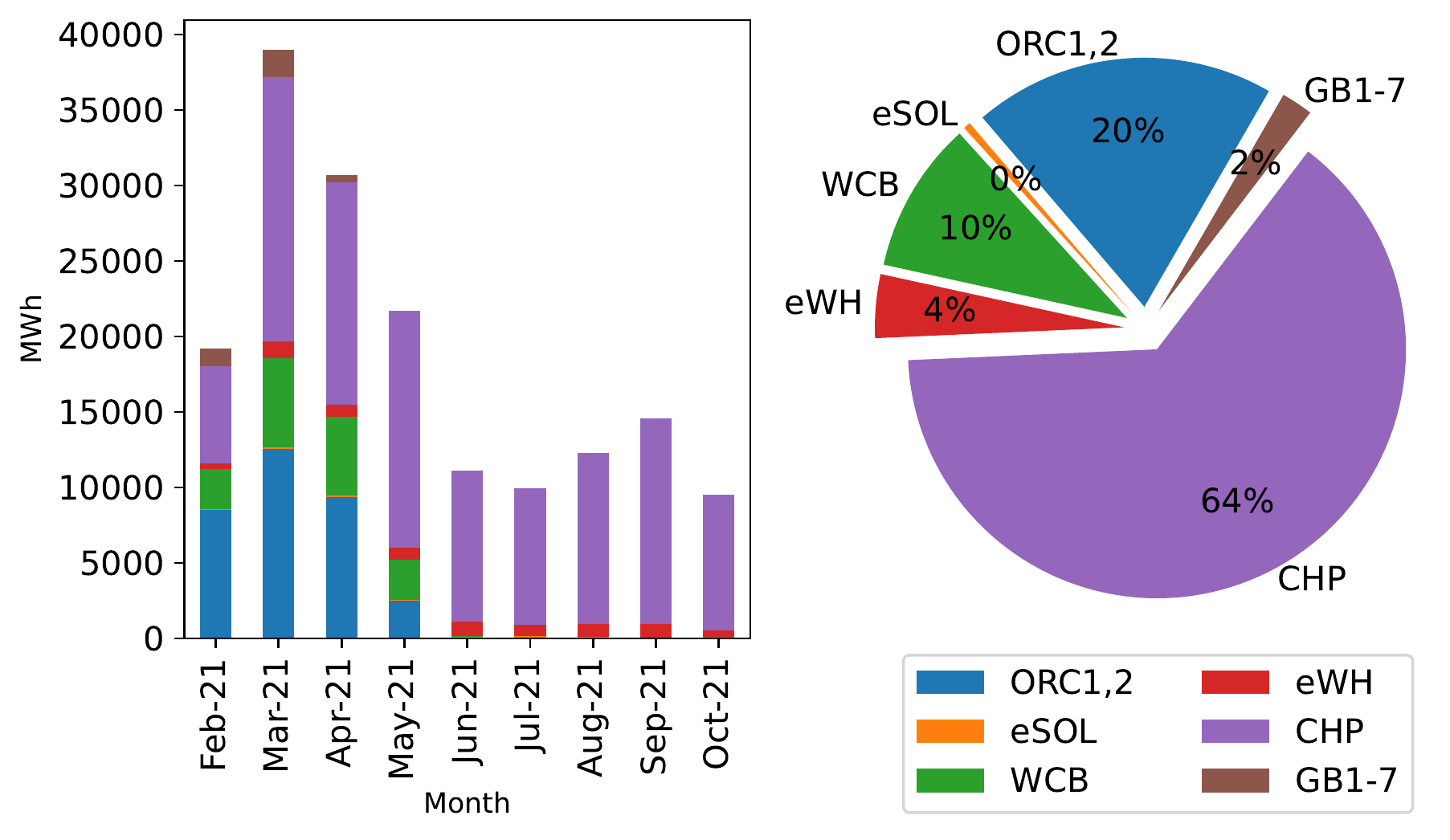}
         \caption{H-02-5808}
         \label{fig:h-total}
     \end{subfigure}
     \begin{subfigure}[b]{0.49\columnwidth}
         \centering
         \includegraphics[width=1\columnwidth]{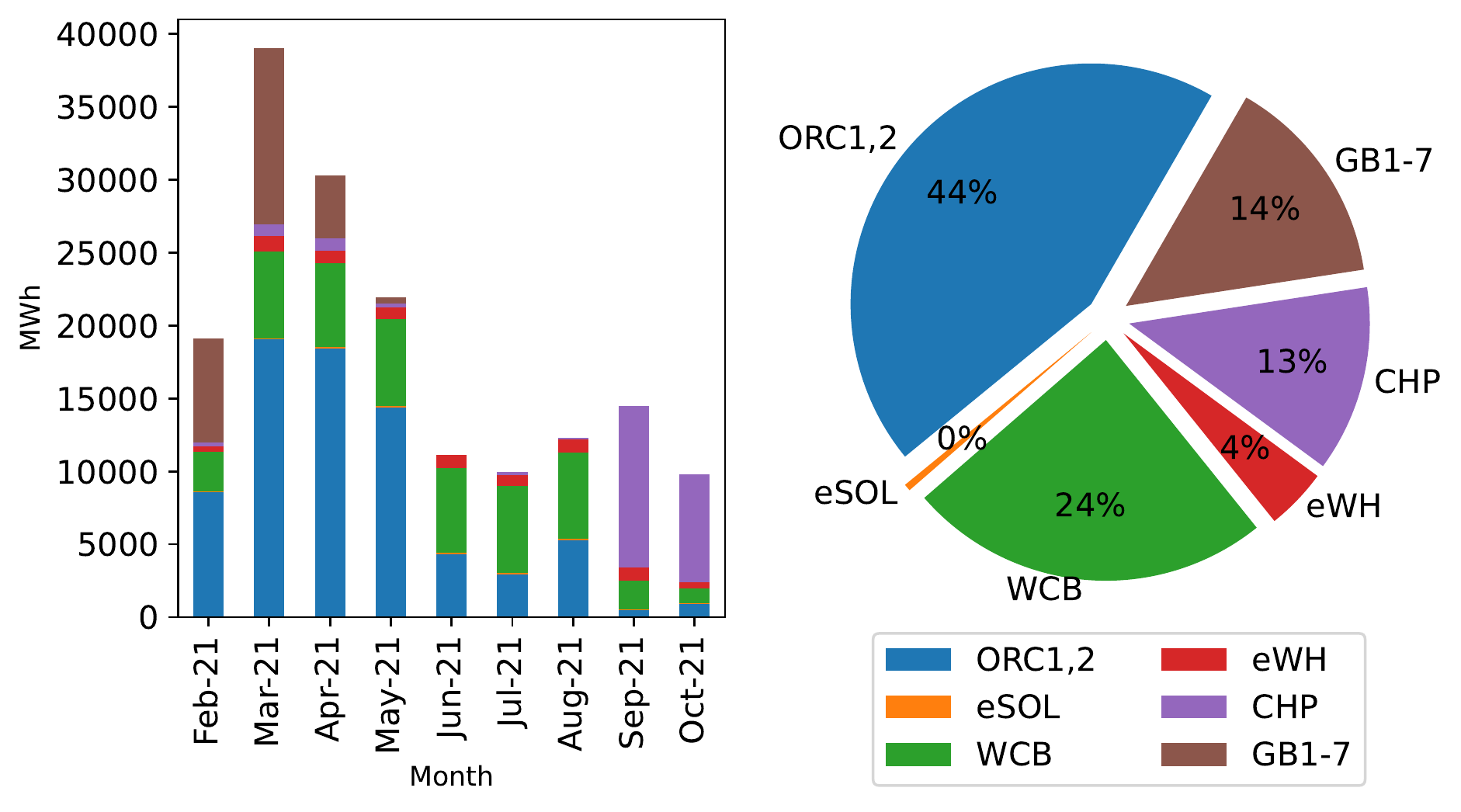}
         \caption{H-02-5808 + 0.04[EUR/kWh] for gas}
         \label{fig:hg-total}
     \end{subfigure}
    \caption{Distribution of heat production over the months and in total.}
    \label{fig:longterm}
    
\end{figure}

The distribution of the heat production among the units (per month and in total) for each DH systems is shown in Figure \ref{fig:longterm}. Objective values, runtimes and RES shares are presented in Table \ref{tab:longterm}. The results put the already observed operational patterns presented in Section \ref{sec:rc} in a yearly context and confirm the observations.

In Brønderslev (Fig. \ref{fig:b-total}), the major share of the heat demand is covered by the ORC unit and WCHPs (which we consider as RES). The CHP units and the electric boiler are mostly used during winter where the heat demand is higher. The share of RES in heat production is more than 78\%, when considering the WCHPs as RES. In Middelfart (Fig. \ref{fig:m-total}), the heat demand during summer can be covered by the wood chip boiler (WCB). In winter, the wood pellet boiler (WPB) is used more extensively. CHP2 only operates occasionally in case of high electricity prices. Here, the share of RES in heat production (WCB+WPB) is even higher with 88\%.

In Hillerød (Fig. \ref{fig:h-total}), the operation strategy shifts between summer and winter. In the colder months, the ORC unit and the wood chip boiler are used for a large share of the heat production in addition to the efficient CHP unit. This coincides with the increasing electricity prices in the second half of 2021. During summer, the system mostly relies on waste heat and the CHP unit, which can operate at favorable power prices (see also the representative week in Section \ref{sec:rc}). Although the CHP unit supports the electricity grid in hours with high prices by with efficient co-generation of heat and power, it relies on natural gas. Therefore, the share of RES (ORC1,2, eSol, WCB, eWH) in Hillerød is low with 34\%. In case market conditions change, the Hillerød system can easily adapt: With increasing gas prices and/or higher emission taxes, the result is expected to change in favor of the ORC unit and WCB. Therefore, we also tested the operation with an increase in natural gas prices of 0.04 EUR/kWh by replacing the cost factor at the energy source for gas. The result is shown in Figure \ref{fig:hg-total}. Now, the production of the ORC unit and wood chip boiler (WCB) is more than doubled. Furthermore, it is less profitable to operate the CHP unit, since the prices on the electricity market are not high enough except in September and October. Thus, high demands are covered also by gas boilers. The share of RES increases from 34\% to 73\%, but also the overall cost increase drastically due to higher fuel costs and lower power market revenues.

The runtimes for the long-term models ranges from 150 seconds to 3760 seconds (slightly more than one hour) which is acceptable for such long planning horizons, since these kind of analyses are only performed occasionally.

\begin{table}[t]
    \centering
    \footnotesize
    \caption{Objective values and runtimes for the long-term analysis}
    \begin{adjustbox}{width=0.8\textwidth}
    \begin{tabular}{lrrr}
    \toprule
        Case & Objective [EUR] & Runtime [s] & Share RES [\%]\\ \midrule
        B-10-6936 & 1795156.73 & 412.55 & 78.61\% \\ 
        M-03-6936 & 710407.53 & 156.94 & 88.18\% \\ 
        H-02-5808 & 522125.61 & 3761.08 & 34.03\% \\ 
        H-02-5808 + 0.04[EUR/kWh] for gas  & 6219317.78 & 293.53 & 73.24\% \\\bottomrule
    \end{tabular}
     \end{adjustbox}
    \label{tab:longterm}
\end{table}

\section{Conclusion}\label{sec:conclusion}
In this paper, we propose a novel model formulation for production optimization in DH systems. The optimization model uses a generic formulation that allows the application in a wide variety of DH systems and allows the integration of scenarios to account for uncertainty. The final model is a two-stage stochastic program based on a network flow model that is build using an easily adaptable generic network structure. The model can be used for different planning problems, such as operational planning under uncertainty, optimization of bids to the day-ahead electricity market and long-term evaluations of operations by exchanging data and non-anticipativity constraints.

The general applicability and performance of the approach is evaluated based on real data from the three Danish DH systems of Brønderslev, Hillerød and Middelfart with different characteristics. The calculation of the VSS alongside out-of-sample evaluations show the benefit of using stochastic programming in operational planning problems including uncertain input data, such as RES production, electricity prices and heat demand. In particular, the modelling of bidding to the day-ahead electricity market profits from the stochastic model. Furthermore, we present the results of applying the model in a rolling horizon setting to realizations of the uncertain data in the three DH systems. The solutions of the operational planning and an additional long-term evaluation on historic data with a deterministic setting allow an analysis of operational strategies, costs and energy mixes. In all cases, model runtimes model are short enough for application in practice.
As final remark we like to add that the diverse data basis for the analysis in this publication uses historical data from years 2019 to 2021, which does not reflect the current (2022) energy price trends and volatility. The specific results in terms costs and energy mixes in the DH systems would be affected, however, the modelling approach remains applicable without restriction. The input data can be replaced with more recent data, which was not available to us at the time of the experiments. 

There are several lines for future research. First, the implementation of the bidding curves in this paper does not model minimum up and down times on market-dependent units. Therefore, an extension with proper block bidding techniques such as proposed in \cite{fleten2007stochastic} and applied in \cite{blockbids} is necessary. Second, the modelling of electricity balancing markets would be a valuable addition. The cases and evaluation in this work focuses on the day-ahead market and uses the balancing market only as an imbalance mechanism and not as an opportunity for trading. Furthermore, the model could be used as a basis for extensive experimental evaluations of different setups in DH systems. 

\section*{Acknowledgments}
 This work is funded by Innovation Fund Denmark through the HEAT 4.0 project (no. 8090-00046B). The authors thank Middelfart Fjernvarme a.m.b.a., Brønderslev Forsyning A/S and Hillerød Forsyning for providing their data.

\section*{Conflict of interest}
\noindent The authors declare that they have no conflict of interest.

\bibliography{bibliography}

\end{document}